\documentclass[12pt,reqno,twoside]{amsart}%

\usepackage{amsmath,amsthm,amscd,amsthm,upref,indentfirst}
\usepackage{amsfonts,mathrsfs}

\usepackage{amssymb,amsbsy,bm}
\usepackage{graphicx}
\usepackage{times}
\usepackage{amsaddr}

\usepackage{soul}
\usepackage{slashbox}

\usepackage[us,long,24hr]{datetime}

\usepackage{mathabx}

\usepackage[usenames,dvipsnames]{color}

\theoremstyle{plain}
\newtheorem{theorem}{Theorem}
\newtheorem{assertion}[theorem]{Assertion}

\newtheorem{proposition}[theorem]{Proposition}

\theoremstyle{definition}
\newtheorem{definition}[theorem]{Definition}

\newtheorem{corollary}[theorem]{Corollary}

\theoremstyle{remark}

\newtheorem{example}[theorem]{Example}
\numberwithin{equation}{section}
\numberwithin{theorem}{section}

\usepackage{exscale}
\usepackage[scaled=0.86]{helvet}
\renewcommand{\mathfrak}[1]{{\textbf{\upshape #1}}}
\renewcommand{\mathbf}{\bm}
\renewcommand{\emph}[1]{\textrm{{\upshape #1}}}
\usepackage[scr=euler,scrscaled=1.15,bb=txof,bbscaled=1.16,cal=boondox,calscaled=1.15]{mathalfa}
\renewcommand{\mathit}[1]{\mathscr #1}
\renewcommand{\mathtt}[1]{\scalebox{1}{\bfseries \texttt{\upshape #1}}}

\usepackage{float}

\usepackage[justification=centering]{caption}

\numberwithin{equation}{section}
\numberwithin{theorem}{section}

\usepackage[normalem]{ulem}
\usepackage[square,comma]{natbib}

\usepackage{mparhack}

\usepackage{keyval}
\usepackage{totcount}
\newtotcounter{citnum} 
\def\oldbibitem{} \let\oldbibitem=\bibitem
\def\bibitem{\stepcounter{citnum}\oldbibitem}

\usepackage[verbose]{backref}
\backrefsetup{verbose=false}
\renewcommand*{\backref}[1]{}
\renewcommand*{\backrefalt}[4]{[{\tiny%
    \ifcase #1 \textsl{Not cited}%
          \or \textsl{Cited on page}~\textcolor{BrickRed}{#2}%
          \else \textsl{Cited on pages}~\textcolor{BrickRed}{#2}%
    \fi%
    }]}

\usepackage{fancyhdr}
\author{\small\scshape S\lowercase{teven} D\lowercase{uplij}}

\address{
Center for Information Technology (WWU IT),
Universit\"at M\"unster,
R\"ontgenstrasse 7-13\\
D-48149 M\"unster,
Deutschland}
\email{\small \sf douplii@uni-muenster.de;
sduplij@gmail.com;
https://ivv5hpp.uni-muenster.de/u/douplii}
\title{\large\bfseries\scshape
P\lowercase{olyadic analogs of direct product
}}

\date{\textit{of completion}
January 19, 2022.
\mbox{}\hskip 1.16em
\textit{Total}:
25
references
}

\renewcommand{\refname}{\textsc{References}}

\let\origsection\section
\renewcommand{\section}[1]{\sectionmark{#1}\origsection{#1}}
\let\origsubsection\subsection
\renewcommand{\subsection}[1]{\subsectionmark{#1}\origsubsection{#1}}

\makeatletter
\renewenvironment{thebibliography}[1]{%
  \@xp\origsection\@xp*\@xp{\refname}%
  \normalfont\footnotesize\labelsep .9em\relax
  \renewcommand\theenumiv{\arabic{enumiv}}\let\p@enumiv\@empty
  \vspace*{-5pt}
  \list{\@biblabel{\theenumiv}}{\settowidth\labelwidth{\@biblabel{#1}}%
    \leftmargin\labelwidth \advance\leftmargin\labelsep
    \usecounter{enumiv}}%
  \sloppy \clubpenalty\@M \widowpenalty\clubpenalty
  \sfcode`\.=\@m
}{%
  \def\@noitemerr{\@latex@warning{Empty `thebibliography' environment}}%
  \endlist
}
\makeatother

\subjclass[2010]{16T25, 17A42, 20B30, 20F36, 20M17, 20N15}
\keywords{direct product, direct power, polyadic semigroup, arity, polyadic ring, polyadic field}

\begin{document}
\mbox{}

\mbox{}
\begin{abstract}

\noindent We propose a generalization of the external direct product concept
to polyadic algebraic structures which introduces novel properties in two
ways: the arity of the product can differ from that of the constituents, and
the elements from different multipliers can be \textquotedblleft
entangled\textquotedblright\ such that the product is no longer componentwise.
The main property which we want to preserve is associativity, which is gained
by using the associativity quiver technique provided earlier. For polyadic
semigroups and groups we introduce two external products: 1) the iterated
direct product which is componentwise, but can have arity different from the
multipliers; 2) the hetero product (power) which is noncomponentwise and
constructed by analogy with the heteromorphism concept introduced earlier. It
is shown in which cases the product of polyadic groups can itself be a
polyadic group. In the same way the external product of polyadic rings and
fields is generalized. The most exotic case is the external product of
polyadic fields, which can be a polyadic field (as opposed to the binary
fields), when all multipliers are zeroless fields. Many illustrative concrete
examples are presented.

\end{abstract}
\maketitle

\thispagestyle{empty}

\mbox{}

\tableofcontents
\newpage

\pagestyle{fancy}

\addtolength{\footskip}{15pt}

\renewcommand{\sectionmark}[1]{%
\markboth{
{ \scshape #1}}{}}

\renewcommand{\subsectionmark}[1]{%
\markright{
\mbox{\;}\\[5pt]
\textmd{#1}}{}}

\fancyhead{}
\fancyhead[EL,OR]{\leftmark}
\fancyhead[ER,OL]{\rightmark}
\fancyfoot[C]{\scshape -- \textcolor{BrickRed}{\thepage} --}

\renewcommand\headrulewidth{0.5pt}
\fancypagestyle {plain1}{ %
\fancyhf{}
\renewcommand {\headrulewidth }{0pt}
\renewcommand {\footrulewidth }{0pt}
}

\fancypagestyle{fancyref}{ %
\fancyhf{} 
\fancyhead[C]{\scshape R\lowercase{eferences} }
\fancyfoot[C]{\scshape -- \textcolor{BrickRed}{\thepage} --}
\renewcommand {\headrulewidth }{0.5pt}
\renewcommand {\footrulewidth }{0pt}
}

\fancypagestyle{emptyf}{
\fancyhead{}
\fancyfoot[C]{\scshape -- \textcolor{BrickRed}{\thepage} --}
\renewcommand{\headrulewidth}{0pt}
}
\mbox{}
\vskip 3.5cm
\thispagestyle{emptyf}

\section{\textsc{Introduction}}

The concept of direct product plays a crucial role for algebraic structures in
the study of their internal constitution and their representation in terms of
better known/simpler structures (see, e.g. \cite{Lang,lambek}). The general
method of the external direct product construction is to take the Cartesian
product of the underlying sets and endow it with the operations from the
algebraic structures under consideration. Usually this is an identical
repetition of the initial multipliers' operations componentwise
\cite{hungerfold}. In the case of polyadic algebraic structures their arity
comes into the game, such that endowing the product with operations becomes
nontrivial in two aspects: the arities of all structures can be different (but
\textquotedblleft quantized\textquotedblright\ and not unique) and the
elements from different multipliers can be \textquotedblleft
entangled\textquotedblright\ making the product not componentwise. The direct
(componentwise) product of $n$-ary groups was considered in
\cite{mic84f,shc2014}. We propose two corresponding polyadic analogs (changing
arity and \textquotedblleft entangling\textquotedblright) of the external
direct product which preserve associativity, and therefore allow us to work
out polyadic semigroups, groups, rings and fields.

The direct product is important, especially because it plays the role of a
product in a corresponding category (see, e.g. \cite{borceaux1,maclane1}). For
instance, the class of all polyadic groups for objects and polyadic group
homomorphisms for morphisms form a category which is well-defined, because it
has the polyadic direct product \cite{mic84a,ian91} as a product.

We then consider polyadic rings and fields in the same way. Since there exist
zeroless polyadic fields \cite{dup2017a}, the well-known statement (see, e.g.
\cite{lambek}) of the absence of binary fields that are a direct product of
fields does not generalize. We construct polyadic fields which are products of
zeroless fields, which can lead to a new category of polyadic fields. The
proposed constructions are accompanied by concrete illustrative examples.

\section{\textsc{Preliminaries}}

We introduce here briefly the usual notation, for details see \cite{dup2018a}.
For a non-empty (underlying) set $G$ the $n$\textit{-tuple} (or
\textit{polyad} \cite{pos}) of elements is denoted by $\left(  g_{1}%
,\ldots,g_{n}\right)  $, $g_{i}\in G$, $i=1,\ldots,n$, and the Cartesian
product is denoted by $G^{\times n}\equiv\overset{n}{\overbrace{G\times
\ldots\times G}}$ and consists of all such $n$-tuples. For all elements equal
to $g\in G$, we denote $n$-tuple (polyad) by a power $\left(  g^{n}\right)  $.
To avoid unneeded indices we denote with one bold letter $\left(
\mathbf{g}\right)  $ a polyad for which the number of elements in the
$n$-tuple is clear from the context, and sometimes we will write $\left(
\mathbf{g}^{\left(  n\right)  }\right)  $. On the Cartesian product $G^{\times
n}$ we define a polyadic (or $n$-ary) operation $\mu^{\left(  n\right)
}:G^{\times n}\rightarrow G$ such that $\mu^{\left(  n\right)  }\left[
\mathbf{g}\right]  \mapsto h$, where $h\in G$. The operations with $n=1,2,3$
are called \textit{unary, binary and ternary}.

Recall the definitions of some algebraic structures and their special elements
(in the notation of \cite{dup2018a}). A (one-set) \textit{polyadic algebraic
structure} $\mathcal{G}$ is a set $G$ closed with respect to polyadic
operations. In the case of one $n$-ary operation $\mu^{\left(  n\right)
}:G^{\times n}\rightarrow G$, it is called \textit{polyadic multiplication}
(or $n$\textit{-ary multiplication}). A one-set $n$\textit{-ary algebraic
structure} $\mathcal{M}^{\left(  n\right)  }=\left\langle G\mid\mu^{\left(
n\right)  }\right\rangle $ or \textit{polyadic magma} ($n$-\textit{ary magma)
}is a set $G$ closed with respect to one $n$-ary operation $\mu^{\left(
n\right)  }$ and without any other additional structure. In the binary case
$\mathcal{M}^{\left(  2\right)  }$ was also called a groupoid by Hausmann and
Ore \cite{hau/ore} (and \cite{cli/pre1}). Since the term \textquotedblleft
groupoid\textquotedblright\ was widely used in category theory for a different
construction, the so-called Brandt groupoid \cite{bra1,bruck}, Bourbaki
\cite{bourbaki98} later introduced the term \textquotedblleft
magma\textquotedblright.

Denote the number of iterating multiplications by $\ell_{\mu}$, and call the
resulting composition an \textit{iterated product} $\left(  \mu^{\left(
n\right)  }\right)  ^{\circ\ell_{\mu}}$, such that%
\begin{equation}
\mu^{\prime\left(  n^{\prime}\right)  }=\left(  \mu^{\left(  n\right)
}\right)  ^{\circ\ell_{\mu}}\overset{def}{=}\overset{\ell_{\mu}}%
{\overbrace{\mu^{\left(  n\right)  }\circ\left(  \mu^{\left(  n\right)  }%
\circ\ldots\left(  \mu^{\left(  n\right)  }\times\operatorname*{id}%
\nolimits^{\times\left(  n-1\right)  }\right)  \ldots\times\operatorname*{id}%
\nolimits^{\times\left(  n-1\right)  }\right)  }}, \label{cg-mn}%
\end{equation}
where the arities are connected by%
\begin{equation}
n^{\prime}=n_{iter}=\ell_{\mu}\left(  n-1\right)  +1, \label{cg-n}%
\end{equation}
which gives the length of a iterated polyad $\left(  \mathbf{g}\right)  $ in
our notation $\left(  \mu^{\left(  n\right)  }\right)  ^{\circ\ell_{\mu}%
}\left[  \mathbf{g}\right]  $.

A \textit{polyadic zero} of a polyadic algebraic structure $\mathcal{G}%
^{\left(  n\right)  }\left\langle G\mid\mu^{\left(  n\right)  }\right\rangle $
is a distinguished element $z\in G$ (and the corresponding $0$-ary operation
$\mu_{z}^{\left(  0\right)  }$) such that for any $\left(  n-1\right)  $-tuple
(polyad) $\mathbf{g}^{\left(  n-1\right)  }\mathbf{\in}G^{\times\left(
n-1\right)  }$ we have%
\begin{equation}
\mu^{\left(  n\right)  }\left[  \mathbf{g}^{\left(  n-1\right)  },z\right]
=z, \label{cg-z}%
\end{equation}
where $z$ can be on any place in the l.h.s. of (\ref{cg-z}). If its place is
not fixed it can be a single zero. As in the binary case, an analog of
positive powers of an element \cite{pos} should coincide with the number of
multiplications $\ell_{\mu}$ in the iteration (\ref{cg-mn}).

A (positive) \textit{polyadic power} of an element is%
\begin{equation}
g^{\left\langle \ell_{\mu}\right\rangle }=\left(  \mu^{\left(  n\right)
}\right)  ^{\circ\ell_{\mu}}\left[  g^{\ell_{\mu}\left(  n-1\right)
+1}\right]  . \label{cg-pp}%
\end{equation}
We define associativity as the invariance of the composition of two $n$-ary
multiplications. An element of a polyadic algebraic structure $g$ is called
$\ell_{\mu}$-\textit{nilpotent} (or simply \textit{nilpotent} for $\ell_{\mu
}=1$), if there exist $\ell_{\mu}$ such that%
\begin{equation}
g^{\left\langle \ell_{\mu}\right\rangle }=z. \label{cg-mz}%
\end{equation}
A \textit{polyadic (}$n$-\textit{ary) identity} (or neutral element) of a
polyadic algebraic structure is a distinguished element $e$ (and the
corresponding $0$-ary operation $\mu_{e}^{\left(  0\right)  }$) such that for
any element $g\in G$ we have%
\begin{equation}
\mu^{\left(  n\right)  }\left[  g,e^{n-1}\right]  =g, \label{cg-e}%
\end{equation}
where $g$ can be on any place in the l.h.s. of (\ref{cg-e}).

In polyadic algebraic structures, there exist \textit{neutral polyads}
$\mathbf{n}\in G^{\times\left(  n-1\right)  }$ satisfying%
\begin{equation}
\mu^{\left(  n\right)  }\left[  g,\mathbf{n}\right]  =g, \label{cg-mng}%
\end{equation}
where $g$ can be on any of $n$ places in the l.h.s. of (\ref{cg-mng}).
Obviously, the sequence of polyadic identities $e^{n-1}$ is a neutral polyad
(\ref{cg-e}).

A one-set polyadic algebraic structure $\left\langle G\mid\mu^{\left(
n\right)  }\right\rangle $ is called \textit{totally associative}, if%
\begin{equation}
\left(  \mu^{\left(  n\right)  }\right)  ^{\circ2}\left[  \mathbf{g}%
,\mathbf{h},\mathbf{u}\right]  =\mu^{\left(  n\right)  }\left[  \mathbf{g}%
,\mu^{\left(  n\right)  }\left[  \mathbf{h}\right]  ,\mathbf{u}\right]
=invariant, \label{cg-ghu}%
\end{equation}
with respect to placement of the internal multiplication $\mu^{\left(
n\right)  }\left[  \mathbf{h}\right]  $ in r.h.s. on any of $n$ places, with a
fixed order of elements in the any fixed polyad of $\left(  2n-1\right)  $
elements $\mathbf{t}^{\left(  2n-1\right)  }=\left(  \mathbf{g},\mathbf{h}%
,\mathbf{u}\right)  \in G^{\times\left(  2n-1\right)  }$.

A \textit{polyadic semigroup} $\mathcal{S}^{\left(  n\right)  }$ is a one-set
$S$ one-operation $\mu^{\left(  n\right)  }$ algebraic structure in which the
$n$-ary multiplication is associative, $\mathcal{S}^{\left(  n\right)
}=\left\langle S\mid\mu^{\left(  n\right)  }\mid\text{associativity
(\ref{cg-ghu})}\right\rangle $. A polyadic algebraic structure $\mathcal{G}%
^{\left(  n\right)  }=\left\langle G\mid\mu^{\left(  n\right)  }\right\rangle
$ is $\sigma$-\textit{commutative}, if $\mu^{\left(  n\right)  }=\mu^{\left(
n\right)  }\circ\sigma$, or%
\begin{equation}
\mu^{\left(  n\right)  }\left[  \mathbf{g}\right]  =\mu^{\left(  n\right)
}\left[  \sigma\circ\mathbf{g}\right]  ,\ \ \ \mathbf{g}\in G^{\times n},
\label{cg-ms}%
\end{equation}
where $\sigma\circ\mathbf{g}=\left(  g_{\sigma\left(  1\right)  }%
,\ldots,g_{\sigma\left(  n\right)  }\right)  $ is a permutated polyad and
$\sigma$ is a fixed element of $S_{n}$, the permutation group on $n$ elements.
If (\ref{cg-ms}) holds for all $\sigma\in S_{n}$, then a polyadic algebraic
structure is \textit{commutative}. A special type of the $\sigma
$-commutativity%
\begin{equation}
\mu^{\left(  n\right)  }\left[  g,\mathbf{t}^{\left(  n-2\right)  },h\right]
=\mu^{\left(  n\right)  }\left[  h,\mathbf{t}^{\left(  n-2\right)  },g\right]
, \label{cg-mth}%
\end{equation}
where $\mathbf{t}^{\left(  n-2\right)  }\in G^{\times\left(  n-2\right)  }$ is
any fixed $\left(  n-2\right)  $-polyad, is called \textit{semicommutativity}.
If an $n$-ary semigroup $\mathcal{S}^{\left(  n\right)  }$ is iterated from a
commutative binary semigroup with identity, then $\mathcal{S}^{\left(
n\right)  }$ is semicommutative. A polyadic algebraic structure is called
(uniquely) $i$-\textit{solvable}, if for all polyads $\mathbf{t}$,
$\mathbf{u}$ and element $h$, one can (uniquely) resolve the equation (with
respect to $h$) for the fundamental operation%
\begin{equation}
\mu^{\left(  n\right)  }\left[  \mathbf{u},h,\mathbf{t}\right]  =g
\label{cg-mug}%
\end{equation}
where $h$ can be on any place, and $\mathbf{u},\mathbf{t}$ are polyads of the
needed length.

A polyadic algebraic structure which is uniquely $i$-solvable for all places
$i=1,\ldots,n$ is called a $n$-\textit{ary }(or \textit{polyadic})\textit{
quasigroup }$\mathcal{Q}^{\left(  n\right)  }=\left\langle Q\mid\mu^{\left(
n\right)  }\mid\text{solvability}\right\rangle $. An associative polyadic
quasigroup is called a $n$-\textit{ary} (or \textit{polyadic})\textit{ group}.
In an $n$-ary group $\mathcal{G}^{\left(  n\right)  }=\left\langle G\mid
\mu^{\left(  n\right)  }\right\rangle $ the only solution of (\ref{cg-mug}) is
called a \textit{querelement} of $g$ and denoted by $\bar{g}$ \cite{dor3},
such that%
\begin{equation}
\mu^{\left(  n\right)  }\left[  \mathbf{h},\bar{g}\right]  =g,\ \ \ g,\bar
{g}\in G, \label{cg-mgg}%
\end{equation}
where $\bar{g}$ can be on any place. Any idempotent $g$ coincides with its
querelement $\bar{g}=g$. The unique solvability relation (\ref{cg-mgg}) in a
$n$-ary group can be treated as a definition of the unary (multiplicative)
\textit{queroperation}%
\begin{equation}
\bar{\mu}^{\left(  1\right)  }\left[  g\right]  =\bar{g}. \label{cg-m1g}%
\end{equation}
We observe from (\ref{cg-mgg}) and (\ref{cg-mng}) that the polyad%
\begin{equation}
\mathbf{n}_{g}=\left(  g^{n-2}\bar{g}\right)  \label{cg-ng}%
\end{equation}
is neutral for any element of a polyadic group, where $\bar{g}$ can be on any
place. If this $i$-th place is important, then we write $\mathbf{n}_{g;i}$. In
a polyadic group the \textit{D\"{o}rnte relations} \cite{dor3}%
\begin{equation}
\mu^{\left(  n\right)  }\left[  g,\mathbf{n}_{h;i}\right]  =\mu^{\left(
n\right)  }\left[  \mathbf{n}_{h;j},g\right]  =g \label{cg-mgnn}%
\end{equation}
hold true for any allowable $i,j$. In the case of a binary group the relations
(\ref{cg-mgnn}) become $g\cdot h\cdot h^{-1}=h\cdot h^{-1}\cdot g=g$.

Using the queroperation (\ref{cg-m1g}) one can give a \textit{diagrammatic
definition} of a polyadic group \cite{gle/gla}: an $n$-\textit{ary group} is a
one-set algebraic structure (universal algebra)%
\begin{equation}
\mathcal{G}^{\left(  n\right)  }=\left\langle G\mid\mu^{\left(  n\right)
},\bar{\mu}^{\left(  1\right)  }\mid\text{associativity (\ref{cg-ghu}),
D\"{o}rnte relations (\ref{cg-mgnn}) }\right\rangle , \label{cg-diam5}%
\end{equation}
where $\mu^{\left(  n\right)  }$ is a $n$-ary associative multiplication and
$\bar{\mu}^{\left(  1\right)  }$ is the queroperation (\ref{cg-m1g}).

\section{\textsc{Polyadic products of semigroups and groups}}

We start from the standard external direct product construction for
semigroups. Then we show that consistent \textquotedblleft
polyadization\textquotedblright\ of the semigroup direct product, which
preserves associativity, can lead to additional properties:%
\begin{align}
&  \text{\textbf{1}) The arities of the polyadic direct product and power can
differ from that of the initial semigroups.}\label{i}\\
&  \text{\textbf{2}) The components of the polyadic power can contain elements
from different multipliers.} \label{ii}%
\end{align}

We use here a vector-like notation for clarity and convenience in passing to
higher arity generalizations. Begin from the direct product of two (binary)
semigroups $\mathcal{G}_{1,2}\equiv\mathcal{G}_{1,2}^{\left(  2\right)
}=\left\langle G_{1,2}\mid\mu_{1,2}^{\left(  2\right)  }\equiv\left(
\cdot_{1,2}\right)  \mid assoc\right\rangle $, where $G_{1,2}$ are underlying
sets, while $\mu_{1,2}^{\left(  2\right)  }$ are multiplications in
$\mathcal{G}_{1,2}$. On the Cartesian product of the underlying sets
$G^{\prime}=G_{1}\times G_{2}$ we define a \textit{direct product}
$\mathcal{G}_{1}\times\mathcal{G}_{2}=\mathcal{G}^{\prime}=\left\langle
G^{\prime}\mid\mathbf{\mu}^{\prime\left(  2\right)  }\equiv\left(
\bullet^{\prime}\right)  \right\rangle $ of the semigroups $\mathcal{G}_{1,2}$
by the componentwise multiplication of the doubles $\mathbf{G}=\left(
\begin{array}
[c]{c}%
g_{1}\\
g_{2}%
\end{array}
\right)  \in G_{1}\times G_{2}$ (being the Kronecker product of doubles in our
notation) , as%
\begin{equation}
\mathbf{G}^{\left(  1\right)  }\bullet^{\prime}\mathbf{G}^{\left(  2\right)
}=\left(
\begin{array}
[c]{c}%
g_{1}\\
g_{2}%
\end{array}
\right)  ^{\left(  1\right)  }\bullet^{\prime}\left(
\begin{array}
[c]{c}%
g_{1}\\
g_{2}%
\end{array}
\right)  ^{\left(  2\right)  }=\left(
\begin{array}
[c]{c}%
g_{1}^{\left(  1\right)  }\cdot_{1}g_{1}^{\left(  2\right)  }\\
g_{2}^{\left(  1\right)  }\cdot_{2}g_{2}^{\left(  2\right)  }%
\end{array}
\right)  , \label{cg-gg}%
\end{equation}
and in the \textquotedblleft polyadic\textquotedblright\ notation%
\begin{equation}
\mathbf{\mu}^{\prime\left(  2\right)  }\left[  \mathbf{G}^{\left(  1\right)
},\mathbf{G}^{\left(  2\right)  }\right]  =\left(
\begin{array}
[c]{c}%
\mu_{1}^{\left(  2\right)  }\left[  g_{1}^{\left(  1\right)  },g_{1}^{\left(
2\right)  }\right] \\[5pt]%
\mu_{2}^{\left(  2\right)  }\left[  g_{2}^{\left(  1\right)  },g_{2}^{\left(
2\right)  }\right]
\end{array}
\right)  . \label{cg-mg12}%
\end{equation}

Obviously, the associativity of $\mathbf{\mu}^{\prime\left(  2\right)  }$
follows immediately from that of $\mu_{1,2}^{\left(  2\right)  }$, because of
the componentwise multiplication in (\ref{cg-mg12}). If $\mathcal{G}_{1,2}$
are groups with the identities $e_{1,2}\in G_{1,2}$, then the identity of the
direct product is the double $\mathbf{E}=\left(
\begin{array}
[c]{c}%
e_{1}\\
e_{2}%
\end{array}
\right)  $, such that $\mathbf{\mu}^{\prime\left(  2\right)  }\left[
\mathbf{E},\mathbf{G}\right]  =\mathbf{\mu}^{\prime\left(  2\right)  }\left[
\mathbf{G},\mathbf{E}\right]  =\mathbf{G}\in\mathcal{G}$.

\subsection{Full polyadic external product}

The \textquotedblleft polyadization\textquotedblright\ of (\ref{cg-mg12}) is straightforward

\begin{definition}
An\textit{ }$n^{\prime}$-ary \textit{full} \textit{direct product} semigroup
$\mathcal{G}^{\prime\left(  n^{\prime}\right)  }=\mathcal{G}_{1}^{\left(
n\right)  }\times\mathcal{G}_{2}^{\left(  n\right)  }$ consists of (two or
$k$) $n$-ary semigroups (of \textsf{the same} arity $n^{\prime}=n$)%
\begin{equation}
\mathbf{\mu}^{\prime\left(  n\right)  }\left[  \mathbf{G}^{\left(  1\right)
},\mathbf{G}^{\left(  2\right)  },\ldots,\mathbf{G}^{\left(  n\right)
}\right]  =\left(
\begin{array}
[c]{c}%
\mu_{1}^{\left(  n\right)  }\left[  g_{1}^{\left(  1\right)  },g_{1}^{\left(
2\right)  },\ldots,g_{1}^{\left(  n\right)  }\right] \\[5pt]%
\mu_{2}^{\left(  n\right)  }\left[  g_{2}^{\left(  1\right)  },g_{2}^{\left(
2\right)  },\ldots,g_{2}^{\left(  n\right)  }\right]
\end{array}
\right)  , \label{cg-mg1}%
\end{equation}
where the (total) polyadic associativity (\ref{cg-ghu}) of $\mu^{\prime\left(
n^{\prime}\right)  }$ is governed by those of the constituent semigroups
$\mathcal{G}_{1}^{\left(  n\right)  }$ and $\mathcal{G}_{2}^{\left(  n\right)
}$ (or $\mathcal{G}_{1}^{\left(  n\right)  }\ldots\mathcal{G}_{k}^{\left(
n\right)  }$) and the componentwise construction (\ref{cg-mg1}).
\end{definition}

If $\mathcal{G}_{1,2}^{\left(  n\right)  }=\left\langle G_{1,2}\mid\mu
_{1,2}^{\left(  n\right)  },\bar{\mu}_{1,2}^{\left(  1\right)  }\right\rangle
$ are $n$-ary groups (where $\bar{\mu}_{1,2}^{\left(  1\right)  }$ are the
unary multiplicative queroperations (\ref{cg-m1g})), then the queroperation
$\mathbf{\bar{\mu}}^{\prime\left(  1\right)  }$ of the full direct product
group $\mathcal{G}^{\prime\left(  n^{\prime}\right)  }=\left\langle G^{\prime
}\equiv G_{1}\times G_{2}\mid\mathbf{\mu}^{\prime\left(  n^{\prime}\right)
},\mathbf{\bar{\mu}}^{\prime\left(  1\right)  }\right\rangle $ ($n^{\prime}%
=n$) is defined componentwise as follows%
\begin{equation}
\mathbf{\bar{G}}\equiv\mathbf{\bar{\mu}}^{\prime\left(  1\right)  }\left[
\mathbf{G}\right]  =\left(
\begin{array}
[c]{c}%
\bar{\mu}_{1}^{\left(  1\right)  }\left[  g_{1}\right] \\
\bar{\mu}_{2}^{\left(  1\right)  }\left[  g_{2}\right]
\end{array}
\right)  ,\ \ \ \text{or \ \ \ }\mathbf{\bar{G}}=\left(
\begin{array}
[c]{c}%
\bar{g}_{1}\\
\bar{g}_{2}%
\end{array}
\right)  , \label{cg-m1gg}%
\end{equation}
which satisfies $\mu^{\prime\left(  n\right)  }\left[  \mathbf{G}%
,\mathbf{G},\ldots,\mathbf{\bar{G}}\right]  =\mathbf{G}$ with $\mathbf{\bar
{G}}$ on any place (cf. (\ref{cg-mgg})).

\begin{definition}
\label{cg-def-der}A\textit{ }full polyadic direct product $\mathcal{G}%
^{\prime\left(  n\right)  }=\mathcal{G}_{1}^{\left(  n\right)  }%
\times\mathcal{G}_{2}^{\left(  n\right)  }$ is called \textit{derived}, if its
constituents $\mathcal{G}_{1}^{\left(  n\right)  }$ and $\mathcal{G}%
_{2}^{\left(  n\right)  }$ are derived, such that the operations $\mu
_{1,2}^{\left(  n\right)  }$ are compositions of the binary operations
$\mu_{1,2}^{\left(  2\right)  }$, correspondingly.
\end{definition}

In the derived case all the operations in (\ref{cg-mg1}) have the form (see
(\ref{cg-mn})--(\ref{cg-n}))%
\begin{equation}
\mu_{1,2}^{\left(  n\right)  }=\left(  \mu_{1,2}^{\left(  2\right)  }\right)
^{\circ\left(  n-1\right)  },\ \ \ \ \ \mathbf{\mu}^{\left(  n\right)
}=\left(  \mathbf{\mu}^{\left(  2\right)  }\right)  ^{\circ\left(  n-1\right)
}. \label{cg-m2}%
\end{equation}

The operations of the derived polyadic semigroup can be written as (cf. the
binary direct product (\ref{cg-gg})--(\ref{cg-mg12}))%
\begin{equation}
\mathbf{\mu}^{\prime\left(  n\right)  }\left[  \mathbf{G}^{\left(  1\right)
},\mathbf{G}^{\left(  2\right)  },\ldots,\mathbf{G}^{\left(  n\right)
}\right]  =\mathbf{G}^{\left(  1\right)  }\bullet^{\prime}\mathbf{G}^{\left(
2\right)  }\bullet^{\prime}\ldots\bullet^{\prime}\mathbf{G}^{\left(  n\right)
}=\left(
\begin{array}
[c]{c}%
g_{1}^{\left(  1\right)  }\cdot_{1}g_{1}^{\left(  2\right)  }\cdot_{1}%
\ldots\cdot_{1}g_{1}^{\left(  n\right)  }\\[5pt]%
g_{2}^{\left(  1\right)  }\cdot_{2}g_{2}^{\left(  2\right)  }\cdot_{2}%
\ldots\cdot_{2}g_{2}^{\left(  n\right)  }%
\end{array}
\right)  .
\end{equation}

We will be more interested in nonderived polyadic analogs of the direct product.

\begin{example}
Let us have two ternary groups: the unitless nonderived group $\mathcal{G}%
_{1}^{\left(  3\right)  }=\left\langle \mathsf{i}\mathbb{R}\mid\mu
_{1}^{\left(  3\right)  }\right\rangle $, where $\mathsf{i}^{2}=-1$, $\mu
_{1}^{\left(  3\right)  }\left[  g_{1}^{\left(  1\right)  },g_{1}^{\left(
2\right)  },g_{1}^{\left(  3\right)  }\right]  =g_{1}^{\left(  1\right)
}g_{1}^{\left(  2\right)  }g_{1}^{\left(  3\right)  }$ is a triple product in
$\mathbb{C}$, the querelement is $\bar{\mu}_{1}^{\left(  1\right)  }\left[
g_{1}\right]  =1/g_{1}$, and $\mathcal{G}_{2}^{\left(  3\right)
}=\left\langle \mathbb{R}\mid\mu_{2}^{\left(  3\right)  }\right\rangle $ with
$\mu_{2}^{\left(  3\right)  }\left[  g_{2}^{\left(  1\right)  },g_{2}^{\left(
2\right)  },g_{2}^{\left(  3\right)  }\right]  =g_{2}^{\left(  1\right)
}\left(  g_{2}^{\left(  2\right)  }\right)  ^{-1}g_{2}^{\left(  3\right)  }$,
the querelement $\bar{\mu}_{2}^{\left(  1\right)  }\left[  g_{2}\right]
=g_{2}$. Then the ternary \textsf{nonderived} full direct product group
becomes $\mathcal{G}^{\prime\left(  3\right)  }=\left\langle \mathsf{i}%
\mathbb{R}\times\mathbb{R}\mid\mathbf{\mu}^{\prime\left(  3\right)
},\mathbf{\bar{\mu}}^{\prime\left(  1\right)  }\right\rangle $, where%
\begin{equation}
\mathbf{\mu}^{\prime\left(  3\right)  }\left[  \mathbf{G}^{\left(  1\right)
},\mathbf{G}^{\left(  2\right)  },\mathbf{G}^{\left(  3\right)  }\right]
=\left(
\begin{array}
[c]{c}%
g_{1}^{\left(  1\right)  }g_{1}^{\left(  2\right)  }g_{1}^{\left(  3\right)
}\\
g_{2}^{\left(  1\right)  }\left(  g_{2}^{\left(  2\right)  }\right)
^{-1}g_{2}^{\left(  3\right)  }%
\end{array}
\right)  ,\ \ \ \ \ \mathbf{\bar{G}}\equiv\mathbf{\bar{\mu}}^{\prime\left(
1\right)  }\left[  \mathbf{G}\right]  =\left(
\begin{array}
[c]{c}%
1/g_{1}\\
g_{2}%
\end{array}
\right)  ,
\end{equation}
which contains \textsf{no identity}, because $\mathcal{G}_{1}^{\left(
3\right)  }$ is unitless and \textsf{nonderived}.
\end{example}

\subsection{Mixed arity iterated product}

In the polyadic case, the following question arises, which cannot even be
stated in the binary case: is it possible to build a version of the
associative direct product such that it can be nonderived and have
\textsf{different arity than the constituent semigroup arities}? The answer is
yes, which leads to two arity changing constructions: componentwise and noncomponentwise.

\begin{enumerate}
\item \textit{Iterated direct product }($\circledast$). In each of the
constituent polyadic semigroups we use the iterating (\ref{cg-mn})
\textsf{componentwise}, but with \textsf{different} numbers of compositions,
because the same number of compositions evidently leads to the iterated
polyadic direct product. In this case the arity of the direct product is
greater than or equal to the arities of the constituents $n^{\prime}\geq
n_{1},n_{2}$.

\item \textit{Hetero product }($\boxtimes$). The polyadic product of $k$
copies of the same $n$-ary semigroup is constructed using the associativity
quiver technique, which \textsf{mixes ("entangles") elements} from different
multipliers, it is \textsf{noncomponentwise} (by analogy with heteromorphisms
in \cite{dup2018a}), and so it can be called a \textit{hetero product} or
\textit{hetero power} (for coinciding multipliers, i.e. constituent polyadic
semigroups or groups). This gives the arity of the hetero product which is
less than or equal to the arities of the equal multipliers $n^{\prime}\leq n$.
\end{enumerate}

In the first componentwise case \textbf{1}), the constituent multiplications
(\ref{cg-mg1}) are composed from the lower arity ones in the componentwise
way, but the initial arities of up and down components can be different (as
opposed to the binary derived case (\ref{cg-m2}))
\begin{equation}
\mu_{1}^{\left(  n\right)  }=\left(  \mu_{1}^{\left(  n_{1}\right)  }\right)
^{\circ\ell_{\mu1}},\ \ \ \ \ \mu_{2}^{\left(  n\right)  }=\left(  \mu
_{2}^{\left(  n_{2}\right)  }\right)  ^{\circ\ell_{\mu2}},\ \ \ \ \ \ 3\leq
n_{1,2}\leq n-1, \label{cg-mnl}%
\end{equation}
where we exclude the limits: the derived case $n_{1,2}=2$ (\ref{cg-m2}) and
the undecomposed case $n_{1,2}=n$ (\ref{cg-mg1}). Since the total size of the
up and down polyads is the same and coincides with the arity of the double
$\mathbf{G}$ multiplication $n^{\prime}$, using (\ref{cg-n}) we obtain the
\textit{arity compatibility} relations%
\begin{equation}
n^{\prime}=\ell_{\mu1}\left(  n_{1}-1\right)  +1=\ell_{\mu2}\left(
n_{2}-1\right)  +1. \label{cg-n1}%
\end{equation}

\begin{definition}
A \textit{mixed arity polyadic iterated direct product }semigroup
$\mathcal{G}^{\prime\left(  n^{\prime}\right)  }=\mathcal{G}_{1}^{\left(
n_{1}\right)  }\circledast\mathcal{G}_{2}^{\left(  n_{2}\right)  }$ consists
of (two) polyadic semigroups $\mathcal{G}_{1}^{\left(  n_{1}\right)  }$ and
$\mathcal{G}_{2}^{\left(  n_{2}\right)  }$ of \textsf{the different} arity
shapes $n_{1}$ and $n_{2}$%
\begin{equation}
\mathbf{\mu}^{\prime\left(  n^{\prime}\right)  }\left[  \mathbf{G}^{\left(
1\right)  },\mathbf{G}^{\left(  2\right)  },\ldots,\mathbf{G}^{\left(
n^{\prime}\right)  }\right]  =\left(
\begin{array}
[c]{c}%
\left(  \mu_{1}^{\left(  n_{1}\right)  }\right)  ^{\circ\ell_{\mu1}}\left[
g_{1}^{\left(  1\right)  },g_{1}^{\left(  2\right)  },\ldots,g_{1}^{\left(
n\right)  }\right] \\[5pt]%
\left(  \mu_{2}^{\left(  n_{2}\right)  }\right)  ^{\circ\ell_{\mu2}}\left[
g_{2}^{\left(  1\right)  },g_{2}^{\left(  2\right)  },\ldots,g_{2}^{\left(
n\right)  }\right]
\end{array}
\right)  , \label{cg-ma}%
\end{equation}
and the arity compatibility relations (\ref{cg-n1}) hold.
\end{definition}

Observe that it is not the case that any two polyadic semigroups can be
composed in the mixed arity polyadic direct product.

\begin{assertion}
If the arity shapes of two polyadic semigroups $\mathcal{G}_{1}^{\left(
n_{1}\right)  }$ and $\mathcal{G}_{2}^{\left(  n_{2}\right)  }$ satisfy the
compatibility condition%
\begin{equation}
a\left(  n_{1}-1\right)  =b\left(  n_{2}-1\right)  =c,\ \ \ \ \ \ \ \ a,b,c\in
\mathbb{N}, \label{cg-ab}%
\end{equation}
then they can form a mixed arity direct product $\mathcal{G}^{\prime\left(
n^{\prime}\right)  }=\mathcal{G}_{1}^{\left(  n_{1}\right)  }\circledast
\mathcal{G}_{2}^{\left(  n_{2}\right)  }$, where $n^{\prime}=c+1$ (\ref{cg-n1}).
\end{assertion}

\begin{example}
In the case of a $4$-ary and $5$-ary semigroups $\mathcal{G}_{1}^{\left(
4\right)  }$ and $\mathcal{G}_{2}^{\left(  5\right)  }$ the direct product
arity of $\mathcal{G}^{\prime\left(  n^{\prime}\right)  }$ is
\textquotedblleft quantized\textquotedblright\ $n^{\prime}=3\ell_{\mu
1}+1=4\ell_{\mu2}+1$, such that%
\begin{align}
n^{\prime}  &  =12k+1=13,25,37,\ldots,\\
\ell_{\mu1}  &  =4k=4,8,12,\ldots,\\
\ell_{\mu2}  &  =3k=3,6,9,\ldots,\ \ \ k\in\mathbb{N},
\end{align}
and only the first mixed arity $13$-ary direct product semigroup
$\mathcal{G}^{\prime\left(  13\right)  }$ is \textsf{nonderived}. If
$\mathcal{G}_{1}^{\left(  4\right)  }$ and $\mathcal{G}_{2}^{\left(  5\right)
}$ are polyadic groups with the queroperations $\bar{\mu}_{1}^{\left(
1\right)  }$ and $\bar{\mu}_{2}^{\left(  1\right)  }$ correspondingly, then
the iterated direct $\mathcal{G}^{\prime\left(  n^{\prime}\right)  }$ is a
polyadic group with the queroperation $\mathbf{\bar{\mu}}^{\prime\left(
1\right)  }$ given in (\ref{cg-m1gg}).
\end{example}

In the same way one can consider the iterated direct product of any number of
polyadic semigroups.

\subsection{Polyadic hetero product}

In the second noncomponentwise case \textbf{2}) we allow multiplying elements
from different components, and therefore we should consider the Cartesian
$k$-power of sets $G^{\prime}=G^{\times k}$ and endow the corresponding
$k$-tuple with a polyadic operation in such a way that associativity of
$\mathcal{G}^{\left(  n\right)  }$ will govern the associativity of the
product $\mathcal{G}^{\prime\left(  n\right)  }$. In other words we construct
a $k$-power of the polyadic semigroup $\mathcal{G}^{\left(  n\right)  }$ such
that the result $\mathcal{G}^{\prime\left(  n^{\prime}\right)  }$ is an
$n^{\prime}$-ary semigroup.

The general structure of the hetero product formally coincides
\textquotedblleft reversely\textquotedblright\ with the main heteromorphism
equation \cite{dup2018a}. The additional parameter which determines the arity
$n^{\prime}$ of the hetero power of the initial $n$-ary semigroup is the
number of intact elements $\ell_{\operatorname*{id}}$. Thus, we arrive at

\begin{definition}
The \textit{hetero }(\textit{\textquotedblleft entangled\textquotedblright%
})\textit{ }$k$\textit{-power} of the $n$-ary semigroup $\mathcal{G}^{\left(
n\right)  }=\left\langle G\mid\mu^{\left(  n\right)  }\right\rangle $ is the
$n^{\prime}$-ary semigroup defined on the $k$-th Cartesian power $G^{\prime
}=G^{\times k}$, such that $\mathcal{G}^{\prime\left(  n^{\prime}\right)
}=\left\langle G^{\prime}\mid\mathbf{\mu}^{\prime\left(  n^{\prime}\right)
}\right\rangle $,%
\begin{equation}
\mathcal{G}^{\prime\left(  n^{\prime}\right)  }=\left(  \mathcal{G}^{\left(
n\right)  }\right)  ^{\boxtimes k}\equiv\overset{k}{\overbrace{\mathcal{G}%
^{\left(  n\right)  }\boxtimes\ldots\boxtimes\mathcal{G}^{\left(  n\right)  }%
}},
\end{equation}
and the $n^{\prime}$-ary multiplication of $k$-tuples $\mathbf{G}^{T}=\left(
g_{1},g_{2},\ldots,g_{k}\right)  \in G^{\times k}$ is given (informally) by%
\begin{equation}
\mathbf{\mu}^{\prime\left(  n^{\prime}\right)  }\left[  \left(
\begin{array}
[c]{c}%
g_{1}\\
\vdots\\
g_{k}%
\end{array}
\right)  ,\ldots,\left(
\begin{array}
[c]{c}%
g_{k\left(  n^{\prime}-1\right)  }\\
\vdots\\
g_{kn^{\prime}}%
\end{array}
\right)  \right]  =\left(
\genfrac{}{}{0pt}{}{\left.
\begin{array}
[c]{c}%
\mu^{\left(  n\right)  }\left[  g_{1},\ldots,g_{n}\right]  ,\\
\vdots\\
\mu^{\left(  n\right)  }\left[  g_{n\left(  \ell_{\mu}-1\right)  }%
,\ldots,g_{n\ell_{\mu}}\right]
\end{array}
\right\}  \ell_{\mu}}{\left.
\begin{array}
[c]{c}%
g_{n\ell_{\mu}+1},\\
\vdots\\
g_{n\ell_{\mu}+\ell_{\operatorname*{id}}}%
\end{array}
\right\}  \ell_{\operatorname*{id}}}%
\right)  ,\ \ \ g_{i}\in G, \label{cg-mmnn}%
\end{equation}
where $\ell_{\operatorname*{id}}$ is the number of \textit{intact elements} in
the r.h.s., and $\ell_{\mu}=k-\ell_{\operatorname*{id}}$ is the number of
multiplications in the resulting $k$-tuple of the direct product. The hetero
power parameters are connected by the \textit{arity changing formula}
\cite{dup2018a}%
\begin{equation}
n^{\prime}=n-\dfrac{n-1}{k}\ell_{\operatorname*{id}}, \label{cg-nnk}%
\end{equation}
with the integer $\dfrac{n-1}{k}\ell_{\operatorname*{id}}\geq1$.
\end{definition}

The concrete placement of elements and multiplications in (\ref{cg-mmnn}) to
obtain the associative $\mathbf{\mu}^{\prime\left(  n^{\prime}\right)  }$ is
governed by the associativity quiver technique \cite{dup2018a}.

There exist important general numerical relations between the parameters of
the twisted direct power $n^{\prime},n,k,\ell_{\operatorname*{id}}$, which
follow from (\ref{cg-mmnn})--(\ref{cg-nnk}). \textsf{First}, there are
non-strict inequalities for them%
\begin{align}
0  &  \leq\ell_{\operatorname*{id}}\leq k-1,\label{cg-l2}\\
\ell_{\mu}  &  \leq k\leq\left(  n-1\right)  \ell_{\mu},\label{cg-lk}\\
2  &  \leq n^{\prime}\leq n. \label{cg-nn}%
\end{align}

\textsf{Second}, the initial and final arities $n$ and $n^{\prime}$ are not
arbitrary, but \textquotedblleft quantized\textquotedblright\ such that the
fraction in (\ref{cg-nnk}) has to be an integer (see \textsc{Table
\ref{cr-T1}}).

\begin{table}[ptb]
\caption{Hetero power \textquotedblleft quantization\textquotedblright.}%
\label{cr-T1}
%
\par
\begin{center}%
\begin{tabular}
[c]{|c|c|c|c|}\hline
$k$ & $\ell_{\mu}$ & $\ell_{\operatorname*{id}}$ & $n/n^{\prime}%
$\\\hline\hline
$2$ & $1$ & $1$ & $%
\begin{array}
[c]{ccccc}%
n= & 3, & 5, & 7, & \ldots\\
n^{\prime}= & 2, & 3, & 4, & \ldots
\end{array}
$\\\hline
$3$ & $1$ & $2$ & $%
\begin{array}
[c]{ccccc}%
n= & 4, & 7, & 10, & \ldots\\
n^{\prime}= & 2, & 3, & 4, & \ldots
\end{array}
$\\\hline
$3$ & $2$ & $1$ & $%
\begin{array}
[c]{ccccc}%
n= & 4, & 7, & 10, & \ldots\\
n^{\prime}= & 3, & 5, & 7, & \ldots
\end{array}
$\\\hline
$4$ & $1$ & $3$ & $%
\begin{array}
[c]{ccccc}%
n= & 5, & 9, & 13, & \ldots\\
n^{\prime}= & 2, & 3, & 4, & \ldots
\end{array}
$\\\hline
$4$ & $2$ & $2$ & $%
\begin{array}
[c]{ccccc}%
n= & 3, & 5, & 7, & \ldots\\
n^{\prime}= & 2, & 3, & 4, & \ldots
\end{array}
$\\\hline
$4$ & $3$ & $1$ & $%
\begin{array}
[c]{ccccc}%
n= & 5, & 9, & 13, & \ldots\\
n^{\prime}= & 4, & 7, & 10, & \ldots
\end{array}
$\\\hline
\end{tabular}
\end{center}
\end{table}

\begin{assertion}
The hetero power is not unique in both directions, if we do not fix the
initial $n$ and final $n^{\prime}$ arities of $\mathcal{G}^{\left(  n\right)
}$ and $\mathcal{G}^{\prime\left(  n^{\prime}\right)  }$.
\end{assertion}

\begin{proof}
This follows from (\ref{cg-mmnn}) and the hetero power \textquotedblleft
quantization\textquotedblright\ \textsc{Table \ref{cr-T1}}.
\end{proof}

The classification of the hetero powers consists of two limiting cases.

\begin{enumerate}
\item \textit{Intactless power}: there are no intact elements $\ell
_{\operatorname*{id}}=0$. The arity of the hetero power reaches its
\textsf{maximum} and coincides with the arity of the initial semigroup
$n^{\prime}=n$ (see \textit{Example} \ref{cg-ex-m3}).

\item \textit{Binary power}: the final semigroup is of lowest arity, i.e.
binary $n^{\prime}=2$. The number of intact elements is (see \textit{Example}
\ref{cg-ex-m2})%
\begin{equation}
\ell_{\operatorname*{id}}=k\frac{n-2}{n-1}.
\end{equation}

\end{enumerate}

\begin{example}
Consider the cubic power of a $4$-ary semigroup $\mathcal{G}^{\prime\left(
3\right)  }=\left(  \mathcal{G}^{\left(  4\right)  }\right)  ^{\boxtimes3}$
with the identity $e$, then the ternary identity triple in $\mathcal{G}%
^{\prime\left(  3\right)  }$ is $\mathbf{E}^{T}=\left(  e,e,e\right)  $, and
therefore this cubic power is a ternary semigroup with identity.
\end{example}

\begin{proposition}
If the initial $n$-ary semigroup $\mathcal{G}^{\left(  n\right)  }$ contains
an identity, then the hetero power $\mathcal{G}^{\prime\left(  n^{\prime
}\right)  }=\left(  \mathcal{G}^{\left(  n\right)  }\right)  ^{\boxtimes k}$
can contain an identity in the \textsf{intactless case} and the Post-like
quiver \cite{dup2018a}. For the binary power $k=2$ only the one-sided identity
is possible.
\end{proposition}

Let us consider some concrete examples.

\begin{example}
\label{cg-ex-m2}Let $\mathcal{G}^{\left(  3\right)  }=\left\langle G\mid
\mu^{\left(  3\right)  }\right\rangle $ be a ternary semigroup, then we can
construct its power $k=2$ (square) of the doubles $\mathbf{G}=\left(
\begin{array}
[c]{c}%
g_{1}\\
g_{2}%
\end{array}
\right)  \in G\times G=G^{\prime}$ in two ways to obtain the associative
hetero power
\begin{equation}
\mathbf{\mu}^{\prime\left(  2\right)  }\left[  \mathbf{G}^{\left(  1\right)
},\mathbf{G}^{\left(  2\right)  }\right]  =\left\{
\begin{array}
[c]{c}%
\left(
\begin{array}
[c]{c}%
\mu^{\left(  3\right)  }\left[  g_{1}^{\left(  1\right)  },g_{2}^{\left(
1\right)  },g_{1}^{\left(  2\right)  }\right] \\
g_{2}^{\left(  2\right)  }%
\end{array}
\right)  ,\\
\left(
\begin{array}
[c]{c}%
\mu^{\left(  3\right)  }\left[  g_{1}^{\left(  1\right)  },g_{2}^{\left(
2\right)  },g_{1}^{\left(  2\right)  }\right] \\
g_{2}^{\left(  1\right)  }%
\end{array}
\right)  ,
\end{array}
\right.  \ \ \ g_{i}^{\left(  j\right)  }\in G. \label{cg-m2g}%
\end{equation}
This means that the Cartesian square can be endowed with the associative
multiplication $\mathbf{\mu}^{\prime\left(  2\right)  }$, and therefore
$\mathcal{G}^{\prime\left(  2\right)  }=\left\langle G^{\prime}\mid
\mathbf{\mu}^{\prime\left(  2\right)  }\right\rangle $ is a binary semigroup
being the hetero product $\mathcal{G}^{\prime\left(  2\right)  }%
=\mathcal{G}^{\left(  3\right)  }\boxtimes\mathcal{G}^{\left(  3\right)  }$.
If $\mathcal{G}^{\left(  3\right)  }$ has a ternary identity $e\in G$, then
$\mathcal{G}^{\prime\left(  2\right)  }$ has only the left (right) identity
$\mathbf{E}=\left(
\begin{array}
[c]{c}%
e\\
e
\end{array}
\right)  \in G^{\prime}$, since $\mathbf{\mu}^{\prime\left(  2\right)
}\left[  \mathbf{E},\mathbf{G}\right]  =\mathbf{G}$ ($\mathbf{\mu}%
^{\prime\left(  2\right)  }\left[  \mathbf{G},\mathbf{E}\right]  =\mathbf{G}%
$), but not the right (left) identity. Thus, $\mathcal{G}^{\prime\left(
2\right)  }$ can be a semigroup only, even if $\mathcal{G}^{\left(  3\right)
}$ is a ternary group.
\end{example}

\begin{example}
\label{cg-ex-m3}Take $\mathcal{G}^{\left(  3\right)  }=\left\langle G\mid
\mu^{\left(  3\right)  }\right\rangle $ a ternary semigroup, then the
multiplication on the double $\mathbf{G}=\left(
\begin{array}
[c]{c}%
g_{1}\\
g_{2}%
\end{array}
\right)  \in G\times G=G^{\prime}$ is ternary and \textsf{noncomponentwise}%
\begin{equation}
\mathbf{\mu}^{\prime\left(  3\right)  }\left[  \mathbf{G}^{\left(  1\right)
},\mathbf{G}^{\left(  2\right)  },\mathbf{G}^{\left(  3\right)  }\right]
=\left(
\begin{array}
[c]{c}%
\mu^{\left(  3\right)  }\left[  g_{1}^{\left(  1\right)  },g_{2}^{\left(
2\right)  },g_{1}^{\left(  3\right)  }\right] \\
\mu^{\left(  3\right)  }\left[  g_{2}^{\left(  1\right)  },g_{1}^{\left(
2\right)  },g_{2}^{\left(  3\right)  }\right]
\end{array}
\right)  ,\ \ \ g_{i}^{\left(  j\right)  }\in G,
\end{equation}
and $\mu^{\prime\left(  3\right)  }$ is associative (and described by the
Post-like associative quiver \cite{dup2018a}), and therefore the cubic hetero
power is the ternary semigroup $\mathcal{G}^{\prime\left(  3\right)
}=\left\langle G\times G\mid\mathbf{\mu}^{\prime\left(  3\right)
}\right\rangle $, such that $\mathcal{G}^{\prime\left(  3\right)
}=\mathcal{G}^{\left(  3\right)  }\boxtimes\mathcal{G}^{\left(  3\right)  }$.
In this case, as opposed to the previous example, the existence of a ternary
identity in $\mathcal{G}^{\left(  3\right)  }$ implies the ternary identity in
the direct cube $\mathcal{G}^{\prime\left(  3\right)  }$ by $\mathbf{E}%
=\left(
\begin{array}
[c]{c}%
e\\
e
\end{array}
\right)  $. If $\mathcal{G}^{\left(  3\right)  }$ is a ternary group with the
unary queroperation $\bar{\mu}^{\left(  1\right)  }$, then the cubic hetero
power $\mathcal{G}^{\prime\left(  3\right)  }$ is also a ternary group of the
special class \cite{dud90}: all querelements coincide (cf. (\ref{cg-m1gg})),
such that $\mathbf{\bar{G}}^{T}=\left(  g_{quer},g_{quer}\right)  $, where
$\bar{\mu}^{\left(  1\right)  }\left[  g\right]  =g_{quer}$, $\forall g\in G$.
This is because in (\ref{cg-mgg}) the querelement can be on any place.
\end{example}

\begin{theorem}
If $\mathcal{G}^{\left(  n\right)  }$ is an $n$-ary group, then the hetero
$k$-power $\mathcal{G}^{\prime\left(  n^{\prime}\right)  }=\left(
\mathcal{G}^{\left(  n\right)  }\right)  ^{\boxtimes k}$ can contain
queroperations in the \textsf{intactless case} only.
\end{theorem}

\begin{corollary}
If the power multiplication (\ref{cg-mmnn}) contains no intact elements
$\ell_{\operatorname*{id}}=0$ and does not change arity $n^{\prime}=n$, a
hetero power can be a polyadic group which has only one querelement.
\end{corollary}

Next we consider more complicated hetero power (\textquotedblleft
entangled\textquotedblright) constructions with and without intact elements,
as well as Post-like and non-Post associative quivers \cite{dup2018a}.

\begin{example}
Let $\mathcal{G}^{\left(  4\right)  }=\left\langle G\mid\mu^{\left(  4\right)
}\right\rangle $ be a $4$-ary semigroup, then we can construct its $4$-ary
associative cubic hetero power $\mathcal{G}^{\prime\left(  4\right)
}=\left\langle G^{\prime}\mid\mathbf{\mu}^{\prime\left(  4\right)
}\right\rangle $ using the Post-like and non-Post associative quivers without
intact elements. Taking in (\ref{cg-mmnn}) $n^{\prime}=n$, $k=3$,
$\ell_{\operatorname*{id}}=0$, we get two possibilities for the multiplication
of the triples $\mathbf{G}^{T}=\left(  g_{1},g_{2},g_{3}\right)  \in G\times
G\times G=G^{\prime}$

\begin{enumerate}
\item \textsf{Post-like associative quiver}. The multiplication of the hetero
cubic power case takes the form%
\begin{equation}
\mathbf{\mu}^{\prime\left(  4\right)  }\left[  \mathbf{G}^{\left(  1\right)
},\mathbf{G}^{\left(  2\right)  },\mathbf{G}^{\left(  3\right)  }%
,\mathbf{G}^{\left(  4\right)  }\right]  =\left(
\begin{array}
[c]{c}%
\mu^{\left(  4\right)  }\left[  g_{1}^{\left(  1\right)  },g_{2}^{\left(
2\right)  },g_{3}^{\left(  3\right)  },g_{1}^{\left(  4\right)  }\right]
\\[5pt]%
\mu^{\left(  4\right)  }\left[  g_{2}^{\left(  1\right)  },g_{3}^{\left(
2\right)  },g_{1}^{\left(  3\right)  },g_{2}^{\left(  4\right)  }\right]
\\[5pt]%
\mu^{\left(  4\right)  }\left[  g_{3}^{\left(  1\right)  },g_{1}^{\left(
2\right)  },g_{2}^{\left(  3\right)  },g_{3}^{\left(  4\right)  }\right]
\end{array}
\right)  ,\ \ \ g_{i}^{\left(  j\right)  }\in G, \label{cg-m4}%
\end{equation}
and it can be shown that $\mathbf{\mu}^{\prime\left(  4\right)  }$ is totally
associative, therefore $\mathcal{G}^{\prime\left(  4\right)  }=\left\langle
G^{\prime}\mid\mathbf{\mu}^{\prime\left(  4\right)  }\right\rangle $ is a
$4$-ary semigroup.

\item \textsf{Non-Post associative quiver}. The multiplication of the hetero
cubic power differs from (\ref{cg-m4})%
\begin{equation}
\mathbf{\mu}^{\prime\left(  4\right)  }\left[  \mathbf{G}^{\left(  1\right)
},\mathbf{G}^{\left(  2\right)  },\mathbf{G}^{\left(  3\right)  }%
,\mathbf{G}^{\left(  4\right)  }\right]  =\left(
\begin{array}
[c]{c}%
\mu^{\left(  4\right)  }\left[  g_{1}^{\left(  1\right)  },g_{3}^{\left(
2\right)  },g_{2}^{\left(  3\right)  },g_{1}^{\left(  4\right)  }\right]
\\[5pt]%
\mu^{\left(  4\right)  }\left[  g_{2}^{\left(  1\right)  },g_{1}^{\left(
2\right)  },g_{3}^{\left(  3\right)  },g_{2}^{\left(  4\right)  }\right]
\\[5pt]%
\mu^{\left(  4\right)  }\left[  g_{3}^{\left(  1\right)  },g_{2}^{\left(
2\right)  },g_{1}^{\left(  3\right)  },g_{3}^{\left(  4\right)  }\right]
\end{array}
\right)  ,\ \ \ g_{i}^{\left(  j\right)  }\in G,
\end{equation}
and it can be shown that $\mathbf{\mu}^{\prime\left(  4\right)  }$ is totally
associative, therefore $\mathcal{G}^{\prime\left(  4\right)  }=\left\langle
G^{\prime}\mid\mathbf{\mu}^{\prime\left(  4\right)  }\right\rangle $ is a
$4$-ary semigroup.
\end{enumerate}

The following is valid for both the above cases. If $\mathcal{G}^{\left(
4\right)  }$ has the $4$-ary identity satisfying%
\begin{equation}
\mu^{\left(  4\right)  }\left[  e,e,e,g\right]  =\mu^{\left(  4\right)
}\left[  e,e,g,e\right]  =\mu^{\left(  4\right)  }\left[  e,g,e,e\right]
=\mu^{\left(  4\right)  }\left[  g,e,e,e\right]  =g,\ \ \ \forall g\in G,
\label{cg-m4e}%
\end{equation}
then the hetero power $\mathcal{G}^{\prime\left(  4\right)  }$ has the $4$-ary
identity%
\begin{equation}
\mathbf{E}=\left(
\begin{array}
[c]{c}%
e\\
e\\
e
\end{array}
\right)  ,\ \ \ e\in G. \label{cg-ee}%
\end{equation}
In the case where $\mathcal{G}^{\left(  3\right)  }$ is a ternary group with
the unary queroperation $\bar{\mu}^{\left(  1\right)  }$, then the cubic
hetero power $\mathcal{G}^{\prime\left(  4\right)  }$ is also a ternary group
with one querelement (cf. \textit{Example} \ref{cg-ex-m3})%
\begin{equation}
\mathbf{\bar{G}}=\overline{\left(
\begin{array}
[c]{c}%
g_{1}\\
g_{2}\\
g_{3}%
\end{array}
\right)  }=\left(
\begin{array}
[c]{c}%
g_{quer}\\
g_{quer}\\
g_{quer}%
\end{array}
\right)  ,\ \ \ g_{quer}\in G,\ \ g_{i}\in G,
\end{equation}
where $g_{quer}=\bar{\mu}^{\left(  1\right)  }\left[  g\right]  $, $\forall
g\in G$.
\end{example}

A more nontrivial example is a cubic hetero power which has different arity to
the initial semigroup.

\begin{example}
Let $\mathcal{G}^{\left(  4\right)  }=\left\langle G\mid\mu^{\left(  4\right)
}\right\rangle $ be a $4$-ary semigroup, then we can construct its ternary
associative cubic hetero power $\mathcal{G}^{\prime\left(  3\right)
}=\left\langle G^{\prime}\mid\mathbf{\mu}^{\prime\left(  3\right)
}\right\rangle $ using the associative quivers with one intact element and two
multiplications \cite{dup2018a}. Taking in (\ref{cg-mmnn}) the parameters
$n^{\prime}=3$, $n=4$, $k=3$, $\ell_{\operatorname*{id}}=1$ (see third line of
\textsc{Table \ref{cr-T1}}), we get for the ternary multiplication
$\mathbf{\mu}^{\prime\left(  3\right)  }$ for the triples $\mathbf{G}%
^{T}=\left(  g_{1},g_{2},g_{3}\right)  \in G\times G\times G=G^{\prime}$ of
the hetero cubic power case the form%
\begin{equation}
\mathbf{\mu}^{\prime\left(  3\right)  }\left[  \mathbf{G}^{\left(  1\right)
},\mathbf{G}^{\left(  2\right)  },\mathbf{G}^{\left(  3\right)  }\right]
=\left(
\begin{array}
[c]{c}%
\mu^{\left(  4\right)  }\left[  g_{1}^{\left(  1\right)  },g_{2}^{\left(
2\right)  },g_{3}^{\left(  3\right)  },g_{1}^{\left(  3\right)  }\right]
\\[5pt]%
\mu^{\left(  4\right)  }\left[  g_{2}^{\left(  1\right)  },g_{3}^{\left(
2\right)  },g_{1}^{\left(  2\right)  },g_{2}^{\left(  3\right)  }\right]
\\[5pt]%
g_{3}^{\left(  1\right)  }%
\end{array}
\right)  ,\ \ \ g_{i}^{\left(  j\right)  }\in G,
\end{equation}
which is totally associative, and therefore the hetero cubic power of $4$-ary
semigroup $\mathcal{G}^{\left(  4\right)  }=\left\langle G\mid\mu^{\left(
4\right)  }\right\rangle $ is a ternary semigroup $\mathcal{G}^{\prime\left(
3\right)  }=\left\langle G^{\prime}\mid\mathbf{\mu}^{\prime\left(  3\right)
}\right\rangle $, such that $\mathcal{G}^{\prime\left(  3\right)  }=\left(
\mathcal{G}^{\left(  4\right)  }\right)  ^{\boxtimes~3}$. If the initial
$4$-ary semigroup $\mathcal{G}^{\left(  4\right)  }$ has the identity
satisfying (\ref{cg-m4e}), then the ternary hetero power $\mathcal{G}%
^{\prime\left(  3\right)  }$ has only the right ternary identity (\ref{cg-ee})
satisfying one relation%
\begin{equation}
\mathbf{\mu}^{\prime\left(  3\right)  }\left[  \mathbf{G},\mathbf{E}%
,\mathbf{E}\right]  =\mathbf{G},\ \ \ \forall\mathbf{G}\in G^{\times3},
\end{equation}
and therefore $\mathcal{G}^{\prime\left(  3\right)  }$ is a ternary semigroup
with a right identity. If $\mathcal{G}^{\left(  4\right)  }$ is a $4$-ary
group with the queroperation $\bar{\mu}^{\left(  1\right)  }$, then the hetero
power $\mathcal{G}^{\prime\left(  3\right)  }$ can only be a ternary semigroup
, because in $\left\langle G^{\prime}\mid\mathbf{\mu}^{\prime\left(  3\right)
}\right\rangle $ we cannot define the standard queroperation \cite{pos}.
\end{example}

\section{\textsc{Polyadic products of rings and fields}}

Now we show that the thorough \textquotedblleft
polyadization\textquotedblright\ of operations can lead to some unexpected new
properties of ring and field external direct products. Recall that in the
binary case the external direct product of fields does not exist at all (see,
e.g., \cite{lambek}). The main new peculiarities of the polyadic case are:

\begin{enumerate}
\item The arity shape of the external product ring and its constituent rings
\textsf{can be different}.

\item The external product of polyadic fields \textsf{can be} a polyadic field.
\end{enumerate}

\subsection{External direct product of binary rings}

First, we recall the ordinary (binary) direct product of rings in notation
which would be convenient to generalize to higher arity structures
\cite{dup2017a}. Let us have two binary rings $\mathcal{R}_{1,2}%
\equiv\mathcal{R}_{1,2}^{\left(  2,2\right)  }=\left\langle R_{1,2}\mid
\nu_{1,2}^{\left(  2\right)  }\equiv\left(  +_{1,2}\right)  ,\mu
_{1,2}^{\left(  2\right)  }\equiv\left(  \cdot_{1,2}\right)  \right\rangle $,
where $R_{1,2}$ are underlying sets, while $\nu_{1,2}^{\left(  2\right)  }$
and $\mu_{1,2}^{\left(  2\right)  }$ are additions and multiplications
(satisfying distributivity) in $\mathcal{R}_{1,2}$, correspondingly. On the
Cartesian product of the underlying sets $R^{\prime}=R_{1}\times R_{2}$ one
defines the \textit{external direct product ring} $\mathcal{R}_{1}%
\times\mathcal{R}_{2}=\mathcal{R}^{\prime}=\left\langle R^{\prime}%
\mid\mathbf{\nu}^{\prime\left(  2\right)  }\equiv\left(  +^{\prime}\right)
,\mathbf{\mu}^{\prime\left(  2\right)  }\equiv\left(  \bullet^{\prime}\right)
\right\rangle $ by the componentwise operations (addition and multiplication)
on the doubles $\mathbf{X}=\left(
\begin{array}
[c]{c}%
x_{1}\\
x_{2}%
\end{array}
\right)  \in R_{1}\times R_{2}$ as follows%
\begin{align}
\mathbf{X}^{\left(  1\right)  }+^{\prime}\mathbf{X}^{\left(  2\right)  }  &
=\left(
\begin{array}
[c]{c}%
x_{1}\\
x_{2}%
\end{array}
\right)  ^{\left(  1\right)  }+^{\prime}\left(
\begin{array}
[c]{c}%
x_{1}\\
x_{2}%
\end{array}
\right)  ^{\left(  2\right)  }\equiv\left(
\begin{array}
[c]{c}%
x_{1}^{\left(  1\right)  }\\
x_{2}^{\left(  1\right)  }%
\end{array}
\right)  +^{\prime}\left(
\begin{array}
[c]{c}%
x_{1}^{\left(  2\right)  }\\
x_{2}^{\left(  2\right)  }%
\end{array}
\right)  =\left(
\begin{array}
[c]{c}%
x_{1}^{\left(  1\right)  }+_{1}x_{1}^{\left(  2\right)  }\\
x_{2}^{\left(  1\right)  }+_{2}x_{2}^{\left(  2\right)  }%
\end{array}
\right)  ,\label{cn-xx12}\\
\mathbf{X}^{\left(  1\right)  }\bullet^{\prime}\mathbf{X}^{\left(  2\right)
}  &  =\left(
\begin{array}
[c]{c}%
x_{1}\\
x_{2}%
\end{array}
\right)  ^{\left(  1\right)  }\bullet^{\prime}\left(
\begin{array}
[c]{c}%
x_{1}\\
x_{2}%
\end{array}
\right)  ^{\left(  2\right)  }=\left(
\begin{array}
[c]{c}%
x_{1}^{\left(  1\right)  }\cdot_{1}x_{1}^{\left(  2\right)  }\\
x_{2}^{\left(  1\right)  }\cdot_{2}x_{2}^{\left(  2\right)  }%
\end{array}
\right)  , \label{cn-xx22}%
\end{align}
or in the polyadic notation (with manifest operations)%
\begin{align}
\mathbf{\nu}^{\prime\left(  2\right)  }\left[  \mathbf{X}^{\left(  1\right)
},\mathbf{X}^{\left(  2\right)  }\right]   &  =\left(
\begin{array}
[c]{c}%
\nu_{1}^{\left(  2\right)  }\left[  x_{1}^{\left(  1\right)  },x_{1}^{\left(
2\right)  }\right] \\[5pt]%
\nu_{2}^{\left(  2\right)  }\left[  x_{2}^{\left(  1\right)  },x_{2}^{\left(
2\right)  }\right]
\end{array}
\right)  ,\label{cn-x1}\\
\mathbf{\mu}^{\prime\left(  2\right)  }\left[  \mathbf{X}^{\left(  1\right)
},\mathbf{X}^{\left(  2\right)  }\right]   &  =\left(
\begin{array}
[c]{c}%
\mu_{1}^{\left(  2\right)  }\left[  x_{1}^{\left(  1\right)  },x_{1}^{\left(
2\right)  }\right] \\[5pt]%
\mu_{2}^{\left(  2\right)  }\left[  x_{2}^{\left(  1\right)  },x_{2}^{\left(
2\right)  }\right]
\end{array}
\right)  . \label{cn-x2}%
\end{align}

The associativity and distributivity of the binary direct product operations
$\mathbf{\nu}^{\prime\left(  2\right)  }$ and $\mathbf{\mu}^{\prime\left(
2\right)  }$ are obviously governed by those of the constituent binary rings
$\mathcal{R}_{1}$ and $\mathcal{R}_{2}$, because of the componentwise
construction on the r.h.s. of (\ref{cn-x1})--(\ref{cn-x2}). In the polyadic
case, the construction of the direct product is not so straightforward and can
have additional unusual peculiarities.

\subsection{Polyadic rings}

Here we recall definitions of polyadic rings \cite{cup,cro2,lee/but} in our
notation \cite{dup2017a,dup2018a}. Consider a polyadic structure $\left\langle
R\mid\mu^{\left(  n\right)  },\nu^{\left(  m\right)  }\right\rangle $ with two
operations on the same set $R$: the $m$-ary addition $\nu^{\left(  m\right)
}:R^{\times m}\rightarrow R$ and the $n$-ary multiplication $\mu^{\left(
n\right)  }:R^{\times n}\rightarrow R$. The \textquotedblleft
interaction\textquotedblright\ between operations can be defined using the
polyadic analog of distributivity.

\begin{definition}
\label{cf-def-dis}The \textit{polyadic distributivity} for $\mu^{\left(
n\right)  }$ and $\nu^{\left(  m\right)  }$ consists of $n$ relations%
\begin{align}
&  \mu^{\left(  n\right)  }\left[  \nu^{\left(  m\right)  }\left[
x_{1},\ldots x_{m}\right]  ,y_{2},y_{3},\ldots y_{n}\right] \nonumber\\
&  =\nu^{\left(  m\right)  }\left[  \mu^{\left(  n\right)  }\left[
x_{1},y_{2},y_{3},\ldots y_{n}\right]  ,\mu^{\left(  n\right)  }\left[
x_{2},y_{2},y_{3},\ldots y_{n}\right]  ,\ldots\mu^{\left(  n\right)  }\left[
x_{m},y_{2},y_{3},\ldots y_{n}\right]  \right] \label{cf-dis1}\\
&  \mu^{\left(  n\right)  }\left[  y_{1},\nu^{\left(  m\right)  }\left[
x_{1},\ldots x_{m}\right]  ,y_{3},\ldots y_{n}\right] \nonumber\\
&  =\nu^{\left(  m\right)  }\left[  \mu^{\left(  n\right)  }\left[
y_{1},x_{1},y_{3},\ldots y_{n}\right]  ,\mu^{\left(  n\right)  }\left[
y_{1},x_{2},y_{3},\ldots y_{n}\right]  ,\ldots\mu^{\left(  n\right)  }\left[
y_{1},x_{m},y_{3},\ldots y_{n}\right]  \right] \label{cf-dis2}\\
&  \vdots\nonumber\\
&  \mu^{\left(  n\right)  }\left[  y_{1},y_{2},\ldots y_{n-1},\nu^{\left(
m\right)  }\left[  x_{1},\ldots x_{m}\right]  \right] \nonumber\\
&  =\nu^{\left(  m\right)  }\left[  \mu^{\left(  n\right)  }\left[
y_{1},y_{2},\ldots y_{n-1},x_{1}\right]  ,\mu^{\left(  n\right)  }\left[
y_{1},y_{2},\ldots y_{n-1},x_{2}\right]  ,\ldots\mu^{\left(  n\right)
}\left[  y_{1},y_{2},\ldots y_{n-1},x_{m}\right]  \right]  , \label{cf-dis3}%
\end{align}
where $x_{i},y_{j}\in R$.
\end{definition}

The operations $\mu^{\left(  n\right)  }$ and $\nu^{\left(  m\right)  }$ are
totally \textit{associative}, if (in the invariance definition
\cite{dup2017a,dup2018a})%
\begin{align}
\nu^{\left(  m\right)  }\left[  \mathbf{u},\nu^{\left(  m\right)  }\left[
\mathbf{v}\right]  ,\mathbf{w}\right]   &  =invariant,\label{cf-as1}\\
\mu^{\left(  n\right)  }\left[  \mathbf{x},\mu^{\left(  n\right)  }\left[
\mathbf{y}\right]  ,\mathbf{t}\right]   &  =invariant, \label{cf-as2}%
\end{align}
where the internal products can be on any place, and $\mathbf{y}\in R^{\times
n}$, $\mathbf{v}\in R^{\times m}$, and the polyads $\mathbf{x}$, $\mathbf{t}$,
$\mathbf{u}$, $\mathbf{w}$ are of the needed lengths. In this way both
algebraic structures $\left\langle R\mid\mu^{\left(  n\right)  }\mid
assoc\right\rangle $ and $\left\langle R\mid\nu^{\left(  m\right)  }\mid
assoc\right\rangle $ are \textit{polyadic semigroups} $\mathcal{S}^{\left(
n\right)  }$ and $\mathcal{S}^{\left(  m\right)  }$.

\begin{definition}
A \textit{polyadic }$\left(  m,n\right)  $-\textit{ring} $\mathcal{R}^{\left(
m,n\right)  }$ is a set $R$ with two operations $\mu^{\left(  n\right)
}:R^{\times n}\rightarrow R$ and $\nu^{\left(  m\right)  }:R^{\times
m}\rightarrow R$, such that:

\begin{enumerate}
\item they are distributive (\ref{cf-dis1})-(\ref{cf-dis3});

\item $\left\langle R\mid\mu^{\left(  n\right)  }\mid assoc\right\rangle $ is
a polyadic semigroup;

\item $\left\langle R\mid\nu^{\left(  m\right)  }\mid
assoc,comm,solv\right\rangle $ is a commutative polyadic group.
\end{enumerate}
\end{definition}

In case the multiplicative semigroup $\left\langle R\mid\mu^{\left(  n\right)
}\mid assoc\right\rangle $ of $\mathcal{R}^{\left(  m,n\right)  }$ is
commutative, $\mu^{\left(  n\right)  }\left[  \mathbf{x}\right]  =\mu^{\left(
n\right)  }\left[  \sigma\circ\mathbf{x}\right]  $, for all $\sigma\in S_{n}$,
then $\mathcal{R}^{\left(  m,n\right)  }$ is called a \textit{commutative
polyadic ring}, and if it contains the identity, then $\mathcal{R}^{\left(
m,n\right)  }$ is a $\left(  m,n\right)  $-\textit{semiring}. A polyadic ring
$\mathcal{R}^{\left(  m,n\right)  }$ is called \textit{derived}, if
$\mathbb{\nu}^{\left(  m\right)  }$ and $\mu^{\left(  n\right)  }$ are
repetitions of the binary addition $\left(  +\right)  $ and multiplication
$\left(  \mathbf{\cdot}\right)  $, while $\left\langle R\mid\left(
\mathbf{+}\right)  \right\rangle $ and $\left\langle R\mid\left(
\mathbf{\cdot}\right)  \right\rangle $ are commutative (binary) group and
semigroup respectively.

\subsection{Full polyadic external direct product of $\left(  m,n\right)
$-rings}

Let us consider the following task: for a given polyadic $\left(  m,n\right)
$-ring $\mathcal{R}^{\prime\left(  m,n\right)  }=\left\langle R^{\prime}%
\mid\mathbf{\nu}^{\prime\left(  m\right)  },\mathbf{\mu}^{\prime\left(
n\right)  }\right\rangle \ $to construct a product of all possible (in arity
shape) constituent rings $\mathcal{R}_{1}^{\left(  m_{1},n_{1}\right)  }$ and
$\mathcal{R}_{2}^{\left(  m_{2},n_{2}\right)  }$. The first-hand
\textquotedblleft polyadization\textquotedblright\ of (\ref{cn-x1}%
)--(\ref{cn-x2}) leads to

\begin{definition}
A \textit{full polyadic direct product ring} $\mathcal{R}^{\prime\left(
m,n\right)  }=\mathcal{R}_{1}^{\left(  m,n\right)  }\times\mathcal{R}%
_{2}^{\left(  m,n\right)  }$ consists of (two) polyadic rings of \textsf{the
same} arity shape, such that%
\begin{align}
\mathbf{\nu}^{\prime\left(  m\right)  }\left[  \mathbf{X}^{\left(  1\right)
},\mathbf{X}^{\left(  2\right)  },\ldots,\mathbf{X}^{\left(  m\right)
}\right]   &  =\left(
\begin{array}
[c]{c}%
\nu_{1}^{\left(  m\right)  }\left[  x_{1}^{\left(  1\right)  },x_{1}^{\left(
2\right)  },\ldots,x_{1}^{\left(  m\right)  }\right] \\[5pt]%
\nu_{2}^{\left(  m\right)  }\left[  x_{2}^{\left(  1\right)  },x_{2}^{\left(
2\right)  },\ldots,x_{2}^{\left(  m\right)  }\right]
\end{array}
\right)  ,\label{cn-y1}\\
\mathbf{\mu}^{\prime\left(  n\right)  }\left[  \mathbf{X}^{\left(  1\right)
},\mathbf{X}^{\left(  2\right)  },\ldots,\mathbf{X}^{\left(  n\right)
}\right]   &  =\left(
\begin{array}
[c]{c}%
\mu_{1}^{\left(  n\right)  }\left[  x_{1}^{\left(  1\right)  },x_{1}^{\left(
2\right)  },\ldots,x_{1}^{\left(  n\right)  }\right] \\[5pt]%
\mu_{2}^{\left(  n\right)  }\left[  x_{2}^{\left(  1\right)  },x_{2}^{\left(
2\right)  },\ldots,x_{2}^{\left(  n\right)  }\right]
\end{array}
\right)  , \label{cn-y2}%
\end{align}
where the polyadic associativity (\ref{cg-ghu}) and polyadic distributivity
(\ref{cf-dis1})--(\ref{cf-dis3}) of the direct product operations
$\nu^{\left(  m\right)  }$ and $\mu^{\left(  n\right)  }$ follow from those of
the constituent rings and the componentwise operations in (\ref{cn-y1}%
)--(\ref{cn-y2}).
\end{definition}

\begin{example}
\label{cn-ex-r23}Consider two $\left(  2,3\right)  $-rings $\mathcal{R}%
_{1}^{\left(  2,3\right)  }=\left\langle \left\{  \mathsf{i}x\right\}  \mid
\nu_{1}^{\left(  2\right)  }=\left(  +\right)  ,\mu_{1}^{\left(  3\right)
}=\left(  \cdot\right)  ,x\in\mathbb{Z},\mathsf{i}^{2}=-1\right\rangle $ and
$\mathcal{R}_{2}^{\left(  2,3\right)  }=\left\langle \left\{  \left(
\begin{array}
[c]{cc}%
0 & a\\
b & 0
\end{array}
\right)  \right\}  \mid\nu_{2}^{\left(  2\right)  }=\left(  +\right)  ,\mu
_{2}^{\left(  3\right)  }=\left(  \cdot\right)  ,a,b\in\mathbb{Z}\right\rangle
$, where $\left(  +\right)  $ and $\left(  \cdot\right)  $ are operations in
$\mathbb{Z}$, then their polyadic direct product on the doubles $\mathbf{X}%
^{T}=\left(  \mathsf{i}x,\left(
\begin{array}
[c]{cc}%
0 & a\\
b & 0
\end{array}
\right)  \right)  \in\left(  \mathsf{i}\mathbb{Z},GL^{\text{adiag}}\left(
2,\mathbb{Z}\right)  \right)  $ is defined by%
\begin{align}
\mathbf{\nu}^{\prime\left(  2\right)  }\left[  \mathbf{X}^{\left(  1\right)
},\mathbf{X}^{\left(  2\right)  }\right]   &  =\left(
\begin{array}
[c]{c}%
\mathsf{i}x^{\left(  1\right)  }+\mathsf{i}x^{\left(  2\right)  }\\
\left(
\begin{array}
[c]{cc}%
0 & a^{\left(  1\right)  }+a^{\left(  2\right)  }\\
b^{\left(  1\right)  }+b^{\left(  2\right)  } & 0
\end{array}
\right)
\end{array}
\right)  ,\label{cn-xx1}\\
\mathbf{\mu}^{\prime\left(  3\right)  }\left[  \mathbf{X}^{\left(  1\right)
},\mathbf{X}^{\left(  2\right)  },\mathbf{X}^{\left(  3\right)  }\right]   &
=\left(
\begin{array}
[c]{c}%
\mathsf{i}x^{\left(  1\right)  }x^{\left(  2\right)  }x^{\left(  3\right)  }\\
\left(
\begin{array}
[c]{cc}%
0 & a^{\left(  1\right)  }b^{\left(  2\right)  }a^{\left(  3\right)  }\\
b^{\left(  1\right)  }a^{\left(  2\right)  }b^{\left(  3\right)  } & 0
\end{array}
\right)
\end{array}
\right)  . \label{cnn-xx2}%
\end{align}
The polyadic associativity and distributivity of the direct product operations
$\mathbf{\nu}^{\prime\left(  2\right)  }$ and $\mathbf{\mu}^{\prime\left(
3\right)  }$ are evident, and therefore $\mathcal{R}^{\left(  2,3\right)
}=\left\langle \left\{  \left(  ix,\left(
\begin{array}
[c]{cc}%
0 & a\\
b & 0
\end{array}
\right)  \right)  \right\}  \mid\mathbf{\nu}^{\prime\left(  2\right)
},\mathbf{\mu}^{\prime\left(  3\right)  }\right\rangle $ is a $\left(
2,3\right)  $-ring $\mathcal{R}^{\left(  2,3\right)  }=\mathcal{R}%
_{1}^{\left(  2,3\right)  }\times\mathcal{R}_{2}^{\left(  2,3\right)  }$.
\end{example}

\begin{definition}
A\textit{ }polyadic direct product $\mathcal{R}^{\left(  m,n\right)  }$ is
called \textit{derived}, if both constituent rings $\mathcal{R}_{1}^{\left(
m,n\right)  }$ and $\mathcal{R}_{2}^{\left(  m,n\right)  }$ are derived, such
that the operations $\nu_{1,2}^{\left(  m\right)  }$ and $\mu_{1,2}^{\left(
n\right)  }$ are compositions of the binary operations $\nu_{1,2}^{\left(
2\right)  }$ and $\mu_{1,2}^{\left(  2\right)  }$, correspondingly.
\end{definition}

So, in the derived case (see (\ref{cg-mn}) all the operations in
(\ref{cn-y1})--(\ref{cn-y2}) have the form (cf. (\ref{cg-m2}))%
\begin{align}
\nu_{1,2}^{\left(  m\right)  }  &  =\left(  \nu_{1,2}^{\left(  2\right)
}\right)  ^{\circ\left(  m-1\right)  },\ \ \ \mu_{1,2}^{\left(  n\right)
}=\left(  \nu_{1,2}^{\left(  2\right)  }\right)  ^{\circ\left(  n-1\right)
},\label{cn-n12}\\
\nu^{\left(  m\right)  }  &  =\left(  \nu^{\left(  2\right)  }\right)
^{\circ\left(  m-1\right)  },\ \ \ \mu^{\left(  n\right)  }=\left(
\nu^{\left(  2\right)  }\right)  ^{\circ\left(  n-1\right)  }. \label{cn-n22}%
\end{align}

Thus, the operations of the derived polyadic ring can be written as (cf. the
binary direct product (\ref{cn-xx12})--(\ref{cn-xx22}))%
\begin{align}
\mathbf{\nu}^{\prime\left(  m\right)  }\left[  \mathbf{X}^{\left(  1\right)
},\mathbf{X}^{\left(  2\right)  },\ldots,\mathbf{X}^{\left(  m\right)
}\right]   &  =\left(
\begin{array}
[c]{c}%
x_{1}^{\left(  1\right)  }+_{1}x_{1}^{\left(  2\right)  }+_{1}\ldots+_{1}%
x_{1}^{\left(  m\right)  }\\[5pt]%
x_{2}^{\left(  1\right)  }+_{2}x_{2}^{\left(  2\right)  }+_{2}\ldots+_{2}%
x_{2}^{\left(  m\right)  }%
\end{array}
\right)  ,\label{cn-d1}\\
\mathbf{\mu}^{\prime\left(  n\right)  }\left[  \mathbf{X}^{\left(  1\right)
},\mathbf{X}^{\left(  2\right)  },\ldots,\mathbf{X}^{\left(  n\right)
}\right]   &  =\left(
\begin{array}
[c]{c}%
x_{1}^{\left(  1\right)  }\cdot_{1}x_{1}^{\left(  2\right)  }\cdot_{1}%
\ldots\cdot_{1}x_{1}^{\left(  n\right)  }\\[5pt]%
x_{2}^{\left(  1\right)  }\cdot_{2}x_{2}^{\left(  2\right)  }\cdot_{2}%
\ldots\cdot_{2}x_{2}^{\left(  n\right)  }%
\end{array}
\right)  , \label{cn-d2}%
\end{align}

The external direct product $\left(  2,3\right)  $-ring $\mathcal{R}^{\left(
2,3\right)  }$ from \textit{Example} \ref{cn-ex-r23} is \textsf{not derived},
because both multiplications $\mu_{1}^{\left(  3\right)  }$ and $\mu
_{2}^{\left(  3\right)  }$ there are nonderived.

\subsection{Mixed arity iterated product of $\left(  m,n\right)  $-rings}

Recall, that some polyadic multiplications can be iterated, i. e. composed
(\ref{cg-mn}) from those of lower arity (\ref{cg-n}), also larger than $2$,
and so being nonderived, in general. The nontrivial \textquotedblleft
polyadization\textquotedblright\ of (\ref{cn-x1})--(\ref{cn-x2}) can arise,
when the composition of the separate (up and down) components in the r.h.s. of
(\ref{cn-y1})--(\ref{cn-y2}) will be different, and therefore the external
product operations on the doubles $\mathbf{X}\in R_{1}\times R_{2}$ cannot be
presented in the iterated form (\ref{cg-mn}).

Let now the constituent operations in (\ref{cn-y1})--(\ref{cn-y2}) be composed
from lower arity corresponding operations, but in different ways for the up
and down components, such that%
\begin{align}
\nu_{1}^{\left(  m\right)  }  &  =\left(  \nu_{1}^{\left(  m_{1}\right)
}\right)  ^{\circ\ell_{\nu1}},\ \ \ \ \ \nu_{2}^{\left(  m\right)  }=\left(
\nu_{2}^{\left(  m_{2}\right)  }\right)  ^{\circ\ell_{\nu2}},\ \ \ \ \ 3\leq
m_{1,2}\leq m-1,\\
\mu_{1}^{\left(  n\right)  }  &  =\left(  \mu_{1}^{\left(  n_{1}\right)
}\right)  ^{\circ\ell_{\mu1}},\ \ \ \ \ \mu_{2}^{\left(  n\right)  }=\left(
\mu_{2}^{\left(  n_{2}\right)  }\right)  ^{\circ\ell_{\mu2}},\ \ \ \ \ \ 3\leq
n_{1,2}\leq n-1,
\end{align}
where we exclude the limits: the derived case $m_{1,2}=n_{1,2}=2$
(\ref{cn-n12})--(\ref{cn-n22}) and the uncomposed case $m_{1,2}=m$,
$n_{1,2}=n$ (\ref{cn-y1})--(\ref{cn-y2}). Since the total size of the up and
down polyads is the same and coincides with the arities of the double addition
$m$ and multiplication $n$, using (\ref{cg-n}) we obtain the \textit{arity
compatibility} relations%
\begin{align}
m  &  =\ell_{\nu1}\left(  m_{1}-1\right)  +1=\ell_{\nu2}\left(  m_{2}%
-1\right)  +1,\label{cn-m1}\\
n  &  =\ell_{\mu1}\left(  n_{1}-1\right)  +1=\ell_{\mu2}\left(  n_{2}%
-1\right)  +1. \label{cn-n1}%
\end{align}

\begin{definition}
\label{cn-def-mix}A \textit{mixed arity polyadic direct product ring}
$\mathcal{R}^{\left(  m,n\right)  }=\mathcal{R}_{1}^{\left(  m_{1}%
,n_{1}\right)  }\circledast\mathcal{R}_{2}^{\left(  m_{2},n_{2}\right)  }$
consists of two polyadic rings of \textsf{the different} arity shape, such
that%
\begin{align}
\mathbf{\nu}^{\prime\left(  m\right)  }\left[  \mathbf{X}^{\left(  1\right)
},\mathbf{X}^{\left(  2\right)  },\ldots,\mathbf{X}^{\left(  m\right)
}\right]   &  =\left(
\begin{array}
[c]{c}%
\left(  \nu_{1}^{\left(  m_{1}\right)  }\right)  ^{\circ\ell_{\nu1}}\left[
x_{1}^{\left(  1\right)  },x_{1}^{\left(  2\right)  },\ldots,x_{1}^{\left(
m\right)  }\right] \\[5pt]%
\left(  \nu_{2}^{\left(  m_{2}\right)  }\right)  ^{\circ\ell_{\nu2}}\left[
x_{2}^{\left(  1\right)  },x_{2}^{\left(  2\right)  },\ldots,x_{2}^{\left(
m\right)  }\right]
\end{array}
\right)  ,\label{cn-ma1}\\
\mathbf{\mu}^{\prime\left(  n\right)  }\left[  \mathbf{X}^{\left(  1\right)
},\mathbf{X}^{\left(  2\right)  },\ldots,\mathbf{X}^{\left(  n\right)
}\right]   &  =\left(
\begin{array}
[c]{c}%
\left(  \mu_{1}^{\left(  n_{1}\right)  }\right)  ^{\circ\ell_{\mu1}}\left[
x_{1}^{\left(  1\right)  },x_{1}^{\left(  2\right)  },\ldots,x_{1}^{\left(
n\right)  }\right] \\[5pt]%
\left(  \mu_{2}^{\left(  n_{2}\right)  }\right)  ^{\circ\ell_{\mu2}}\left[
x_{2}^{\left(  1\right)  },x_{2}^{\left(  2\right)  },\ldots,x_{2}^{\left(
n\right)  }\right]
\end{array}
\right)  , \label{cn-ma2}%
\end{align}
and the arity compatibility relations (\ref{cn-m1})--(\ref{cn-n1}) hold valid.
\end{definition}

Thus, two polyadic rings cannot always be composed in the mixed arity polyadic
direct product.

\begin{assertion}
If the arity shapes of two polyadic rings $\mathcal{R}_{1}^{\left(
m_{1},n_{1}\right)  }$ and $\mathcal{R}_{2}^{\left(  m_{2},n_{2}\right)  }$
satisfies the compatibility conditions%
\begin{align}
a\left(  m_{1}-1\right)   &  =b\left(  m_{2}-1\right)  ,\label{cn-ab1}\\
c\left(  n_{1}-1\right)   &  =d\left(  n_{2}-1\right)
,\ \ \ \ \ \ \ \ a,b,c,d\in\mathbb{N}, \label{cn-ab2}%
\end{align}
then they can form a mixed arity direct product.
\end{assertion}

The limiting cases, undecomposed (\ref{cn-y1})--(\ref{cn-y2}) and derived
(\ref{cn-d1})--(\ref{cn-d2}), satisfy the compatibility conditions
(\ref{cn-ab1})--(\ref{cn-ab2}) as well.

\begin{example}
Let us consider two (nonderived) polyadic rings $\mathcal{R}_{1}^{\left(
9,3\right)  }=\left\langle \left\{  8l+7\right\}  \mid\nu_{1}^{\left(
9\right)  },\mu_{1}^{\left(  3\right)  },l\in\mathbb{Z}\right\rangle $ and
$\mathcal{R}_{2}^{\left(  5,5\right)  }=\left\langle \left\{  M\right\}
\mid\nu_{2}^{\left(  5\right)  },\mu_{2}^{\left(  5\right)  }\right\rangle $,
where%
\begin{equation}
M=\left(
\begin{array}
[c]{cccc}%
0 & 4k_{1}+3 & 0 & 0\\
0 & 0 & 4k_{2}+3 & 0\\
0 & 0 & 0 & 4k_{3}+3\\
4k_{4}+3 & 0 & 0 & 0
\end{array}
\right)  ,\ \ \ k_{i}\in\mathbb{Z}, \label{cn-m44}%
\end{equation}
and $\nu_{2}^{\left(  5\right)  }$ and $\mu_{2}^{\left(  5\right)  }$ are the
ordinary sum and product of $5$ matrices. Using (\ref{cn-m1})--(\ref{cn-n1})
we obtain $m=9$, $n=5$, if we choose the smallest \textquotedblleft numbers of
multiplications\textquotedblright\ $\ell_{\nu1}=1$, $\ell_{\nu2}=2$,
$\ell_{\mu1}=2$, $\ell_{\mu2}=1$, and therefore the mixed arity direct product
nonderived $\left(  9,5\right)  $-ring becomes%
\begin{equation}
\mathcal{R}^{\left(  9,5\right)  }=\left\langle \left\{  \mathbf{X}\right\}
\mid\mathbf{\nu}^{\prime\left(  9\right)  },\mathbf{\mu}^{\prime\left(
5\right)  }\right\rangle , \label{cn-r95}%
\end{equation}
where the doubles are $\mathbf{X=}\left(
\begin{array}
[c]{c}%
8l+7\\
M
\end{array}
\right)  $ and the nonderived direct product operations are%
\begin{equation}
\mathbf{\nu}^{\prime\left(  9\right)  }\left[  \mathbf{X}^{\left(  1\right)
},\mathbf{X}^{\left(  2\right)  },\mathbf{\ldots},\mathbf{X}^{\left(
9\right)  }\right]  =\left(
\begin{array}
[c]{c}%
8\left(  l^{\left(  1\right)  }+l^{\left(  2\right)  }+l^{\left(  3\right)
}+l^{\left(  4\right)  }+l^{\left(  5\right)  }+l^{\left(  6\right)
}+l^{\left(  7\right)  }+l^{\left(  8\right)  }+l^{\left(  9\right)
}+7\right)  +7\\
\left(
\begin{array}
[c]{cccc}%
0 & 4K_{1}+3 & 0 & 0\\
0 & 0 & 4K_{2}+3 & 0\\
0 & 0 & 0 & 4K_{3}+3\\
4K_{4}+3 & 0 & 0 & 0
\end{array}
\right)
\end{array}
\right)  ,
\end{equation}%
\begin{equation}
\mathbf{\mu}^{\prime\left(  5\right)  }\left[  \mathbf{X}^{\left(  1\right)
},\mathbf{X}^{\left(  2\right)  },\mathbf{X}^{\left(  3\right)  }%
,\mathbf{X}^{\left(  4\right)  },\mathbf{X}^{\left(  5\right)  }\right]
=\left(
\begin{array}
[c]{c}%
\left(  8l_{\mu}+7\right) \\
\left(
\begin{array}
[c]{cccc}%
0 & 4K_{\mu,1}+3 & 0 & 0\\
0 & 0 & 4K_{\mu,2}+3 & 0\\
0 & 0 & 0 & 4K_{\mu,3}+3\\
4K_{\mu,4}+3 & 0 & 0 & 0
\end{array}
\right)
\end{array}
\right)  ,
\end{equation}
where in the first line $K_{i}=k_{i}^{\left(  1\right)  }+k_{i}^{\left(
2\right)  }+k_{i}^{\left(  3\right)  }+k_{i}^{\left(  4\right)  }%
+k_{i}^{\left(  5\right)  }+k_{i}^{\left(  6\right)  }+k_{i}^{\left(
7\right)  }+k_{i}^{\left(  8\right)  }+k_{i}^{\left(  9\right)  }%
+6\in\mathbb{Z}$, $l_{\mu}\in\mathbb{Z}$ is a cumbersome integer function of
$l^{\left(  j\right)  }\in\mathbb{Z}$, $j=1,\ldots,9$, and in the second line
$K_{\mu,i}\in\mathbb{Z}$ are cumbersome integer functions of $k_{i}^{\left(
s\right)  }$, $i=1,\ldots,4$, $s=1,\ldots,5$. Therefore the polyadic ring
(\ref{cn-r95}) is the nonderived \textsf{mixed }arity polyadic external
product $\mathcal{R}^{\left(  9,5\right)  }=\mathcal{R}_{1}^{\left(
9,3\right)  }\circledast\mathcal{R}_{2}^{\left(  5,5\right)  }$ (see
\textbf{Definition} \ref{cn-def-mix}).
\end{example}

\begin{theorem}
The category of \textsf{polyadic rings} $\mathtt{PolRing}$ can exist (having
the class of all polyadic rings for objects and ring homomorphisms for
morphisms) and can be well-defined, because it has a product as the polyadic
external product of rings.
\end{theorem}

In the same way one can construct the iterated full and mixed arity products
of any number $k$ of polyadic rings, just by passing from the doubles
$\mathbf{X}$ to $k$-tuples $\mathbf{X}_{k}^{T}=\left(  x_{1},\ldots
,x_{k}\right)  $.

\subsection{Polyadic hetero product of $\left(  m,n\right)  $-fields}

The most crucial difference between the binary direct products and the
polyadic ones arises for fields, because a direct product two binary fields is
not a field \cite{lambek}. The reason lies in the fact that each binary field
$\mathcal{F}^{\left(  2,2\right)  }$ necessarily contains $0$ and $1$, by
definition. As follows from (\ref{cn-xx22}), a binary direct product contains
nonzero idempotent doubles $\left(
\begin{array}
[c]{c}%
1\\
0
\end{array}
\right)  $ and $\left(
\begin{array}
[c]{c}%
0\\
1
\end{array}
\right)  $ which are noninvertible, and therefore the external direct product
of fields $\mathcal{F}_{1}^{\left(  2,2\right)  }\times\mathcal{F}%
_{2}^{\left(  2,2\right)  }$ can never be a field. As opposite, polyadic
fields (see \textbf{Definition} \ref{cn-def-fmn}) can be zeroless (we denote
them by hat $\widehat{\mathcal{F}}$), and the above arguments do not hold
valid for them.

Recall definitions of $\left(  m,n\right)  $-fields (see
\cite{lee/but,ian/pop97}). Denote $R^{\ast}=R\setminus\left\{  z\right\}  $,
if the zero $z$ exists (\ref{cg-z}). Observe that (in distinction to binary
rings) $\left\langle R^{\ast}\mid\mu^{\left(  n\right)  }\mid
assoc\right\rangle $ is not a polyadic group, in general. \label{cn-def-dmn}If
$\left\langle R^{\ast}\mid\mu^{\left(  n\right)  }\right\rangle $ is the
$n$-ary group, then $\mathcal{R}^{\left(  m,n\right)  }$ is called a $\left(
m,n\right)  $-\textit{division} \textit{ring} $\mathcal{D}^{\left(
m,n\right)  }$.

\begin{definition}
\label{cn-def-fmn}A (totally) commutative $\left(  m,n\right)  $-division ring
$\mathcal{R}^{\left(  m,n\right)  }$ is called a $\left(  m,n\right)
$-\textit{field }$\mathcal{F}^{\left(  m,n\right)  }$.
\end{definition}

In $n$-ary groups there exists an \textquotedblleft
intermediate\textquotedblright\ commutativity, so called semicommutativity
(\ref{cg-mth}).

\begin{definition}
\label{cn-def-fmn1}A semicommutative $\left(  m,n\right)  $-division ring
$\mathcal{R}^{\left(  m,n\right)  }$ is called a \textit{semicommutative}
$\left(  m,n\right)  $-\textit{field }$\mathcal{F}^{\left(  m,n\right)  }$.
\end{definition}

The definition of a polyadic field can be done in a diagrammatic form,
analogous to (\ref{cg-diam5}). We introduce the \textit{double D\"{o}rnte
relations}: for $n$-ary multiplication $\mu^{\left(  n\right)  }$
(\ref{cg-mgnn}) and for $m$-ary addition $\nu^{\left(  m\right)  }$, as
follows%
\begin{equation}
\nu^{\left(  m\right)  }\left[  \mathbf{m}_{y},x\right]  =x, \label{cn-num}%
\end{equation}
where the (additive) neutral sequence is $\mathbf{m}_{y}=\left(
y^{m-2},\tilde{y}\right)  $, and $\tilde{y}$ is the additive querelement for
$y\in R$ (see (\ref{cg-ng})). As distinct from (\ref{cg-mgnn}) we have only
one (additive) D\"{o}rnte relation (\ref{cn-num}) and one diagram from
(\ref{cg-diam5}) only, because of commutativity of $\nu^{\left(  m\right)  }$.

By analogy with the multiplicative queroperation $\bar{\mu}^{\left(  1\right)
}$ (\ref{cg-m1g}), introduce the \textit{additive unary queroperation} by%
\begin{equation}
\tilde{\nu}^{\left(  1\right)  }\left(  x\right)  =\tilde{x},\ \ \ \ \forall
x\in R, \label{cn-v1}%
\end{equation}
where $\tilde{x}$ is the additive querelement (\ref{cg-m1g}). Thus, we have

\begin{definition}
[\textsf{Diagrammatic definition of }$\left(  m,n\right)  $\textsf{-field}]A
\textit{(polyadic) }$\left(  m,n\right)  $\textit{-field} is a one-set
algebraic structure with 4 operations and 3 relations
\begin{equation}
\left\langle R\mid\nu^{\left(  m\right)  },\tilde{\nu}^{\left(  1\right)
},\mu^{\left(  n\right)  },\bar{\mu}^{\left(  1\right)  }\mid
\text{associativity, distributivity, double D\"{o}rnte relations}\right\rangle
,
\end{equation}
where $\nu^{\left(  m\right)  }$ and $\mu^{\left(  n\right)  }$ are
commutative associative $m$-ary addition and $n$-ary associative
multiplication connected by polyadic distributivity (\ref{cf-dis1}%
)--(\ref{cf-dis3}), $\tilde{\nu}^{\left(  1\right)  }$ and $\bar{\mu}^{\left(
1\right)  }$ are unary additive queroperation\textit{ }(\ref{cn-v1}) and
multiplicative queroperation (\ref{cg-m1g}).
\end{definition}

There is no initial relation between $\tilde{\nu}^{\left(  1\right)  }$ and
$\bar{\mu}^{\left(  1\right)  }$, nevertheless a possible their
\textquotedblleft interaction\textquotedblright\ can lead to further thorough
classification of polyadic fields.

\begin{definition}
\label{cn-def-qsym}A polyadic field $\mathcal{F}^{\left(  m,n\right)  }$ is
called \textit{quer-symmetric}, if its unary queroperations commute%
\begin{align}
\tilde{\nu}^{\left(  1\right)  }\circ\bar{\mu}^{\left(  1\right)  }  &
=\bar{\mu}^{\left(  1\right)  }\circ\tilde{\nu}^{\left(  1\right)
},\label{cn-n1m}\\
\widetilde{\overline{x}}  &  =\overline{\widetilde{x}},\ \ \ \forall x\in R,
\label{cn-n1m1}%
\end{align}
in other case $\mathcal{F}^{\left(  m,n\right)  }$ is called
\textit{quer-nonsymmetric}.
\end{definition}

\begin{example}
\label{cf-ex-tern}Consider the nonunital zeroless (denoted by hat
$\widehat{\mathcal{F}}$) polyadic field $\widehat{\mathcal{F}}^{\left(
3,3\right)  }=\left\langle \left\{  ia/b\right\}  \mid\nu^{\left(  3\right)
},\mu^{\left(  3\right)  }\right\rangle $, $\mathsf{i}^{2}=-1$, $a,b\in
\mathbb{Z}^{odd}$. The ternary addition $\nu^{\left(  3\right)  }\left[
x,y,t\right]  =x+y+t$ and the ternary multiplication $\mu^{\left(  3\right)
}\left[  x,y,t\right]  =xyt$ are nonderived, ternary associative and
distributive (operations are in $\mathbb{C}$). For each $x=ia/b$
($a,b\in\mathbb{Z}^{odd}$) the additive querelement is $\tilde{x}=-ia/b$, and
the multiplicative querelement is $\bar{x}=-ib/a$ (see (\ref{cg-mgg})).
Therefore, both $\left\langle \left\{  ia/b\right\}  \mid\mu^{\left(
3\right)  }\right\rangle $ and $\left\langle \left\{  ia/b\right\}  \mid
\nu^{\left(  3\right)  }\right\rangle $ are ternary groups, but they both
contain \textsf{no neutral elements} (no unit, no zero).The nonunital zeroless
$\left(  3,3\right)  $-field $\widehat{\mathcal{F}}^{\left(  3,3\right)  }$ is
\textsf{quer-symmetric}, because (see (\ref{cn-n1m1}))
\begin{equation}
\widetilde{\overline{x}}=\overline{\widetilde{x}}=i\frac{b}{a}.
\end{equation}

\end{example}

To find quer-nonsymmetric polyadic fields is not a simple task.

\begin{example}
Consider the set of real $4\times4$ matrices over the fractions $\frac
{4k+3}{4l+3}$, $k,l\in\mathbb{Z}$, of the form%
\begin{equation}
M=\left(
\begin{array}
[c]{cccc}%
0 & \dfrac{4k_{1}+3}{4l_{1}+3} & 0 & 0\\
0 & 0 & \dfrac{4k_{2}+3}{4l_{2}+3} & 0\\
0 & 0 & 0 & \dfrac{4k_{3}+3}{4l_{3}+3}\\
\dfrac{4k_{4}+3}{4l_{4}+3} & 0 & 0 & 0
\end{array}
\right)  ,\ \ \ k_{i},l_{i}\in\mathbb{Z}. \label{cn-m4}%
\end{equation}
The set $\left\{  M\right\}  $ is closed with respect to the ordinary addition
of $m\geq5$ matrices, because the sum of feweer of the fractions $\frac
{4k+3}{4l+3}$ does not give a fraction of the same form \cite{dup2017a}, and
with respect to the ordinary multiplication of $n\geq5$ matrices, since the
product of fewer matrices (\ref{cn-m4}) does not have the same shape
\cite{dup2021b}. The polyadic associativity and polyadic distributivity follow
from the binary ones of the ordinary matrices over $\mathbb{R}$, and the
product of $5$ matrices is semicommutative (see \ref{cg-mth}). Taking the
minimal values $m=5$, $n=5$, we define the semicommutative zeroless $\left(
5,5\right)  $-field (see (\ref{cn-def-fmn1}))%
\begin{equation}
\mathcal{F}_{M}^{\left(  5,5\right)  }=\left\langle \left\{  M\right\}
\mid\nu^{\left(  5\right)  },\mu^{\left(  5\right)  },\tilde{\nu}^{\left(
1\right)  },\bar{\mu}^{\left(  1\right)  }\right\rangle , \label{cn-f55}%
\end{equation}
where $\nu^{\left(  5\right)  }$ and $\mu^{\left(  5\right)  }$ are the
ordinary sum and product of $5$ matrices, while $\tilde{\nu}^{\left(
1\right)  }$ and $\bar{\mu}^{\left(  1\right)  }$ are additive and
multiplicative queroperations%
\begin{equation}
\tilde{\nu}^{\left(  1\right)  }\left[  M\right]  \equiv\tilde{M}%
=-3M,\ \ \bar{\mu}^{\left(  1\right)  }\left[  M\right]  \equiv\bar{M}%
=\frac{4l_{1}+3}{4k_{1}+3}\frac{4l_{2}+3}{4k_{2}+3}\frac{4l_{3}+3}{4k_{3}%
+3}\frac{4l_{4}+3}{4k_{4}+3}M. \label{cn-vm}%
\end{equation}
The division ring $\mathcal{D}_{M}^{\left(  5,5\right)  }$ is zeroless,
because the fraction $\frac{4k+3}{4l+3}$, is never zero for $k,l\in\mathbb{Z}%
$, and it is unital with the unit%
\begin{equation}
M_{e}=\left(
\begin{array}
[c]{cccc}%
0 & 1 & 0 & 0\\
0 & 0 & 1 & 0\\
0 & 0 & 0 & 1\\
1 & 0 & 0 & 0
\end{array}
\right)  .
\end{equation}
Using (\ref{cn-m4}) and (\ref{cn-vm}), we obtain%
\begin{align}
\tilde{\nu}^{\left(  1\right)  }\left[  \bar{\mu}^{\left(  1\right)  }\left[
M\right]  \right]   &  =-3\frac{4l_{1}+3}{4k_{1}+3}\frac{4l_{2}+3}{4k_{2}%
+3}\frac{4l_{3}+3}{4k_{3}+3}\frac{4l_{4}+3}{4k_{4}+3}M,\\
\bar{\mu}^{\left(  1\right)  }\left[  \tilde{\nu}^{\left(  1\right)  }\left[
M\right]  \right]   &  =-\frac{1}{27}\frac{4l_{1}+3}{4k_{1}+3}\frac{4l_{2}%
+3}{4k_{2}+3}\frac{4l_{3}+3}{4k_{3}+3}\frac{4l_{4}+3}{4k_{4}+3}M,
\end{align}
or%
\begin{equation}
\widetilde{\overline{M}}=81\overline{\widetilde{M}}, \label{cf-m81}%
\end{equation}
and therefore the additive and multiplicative queroperations do not commute
independently of the field parameters. Thus, the matrix $\left(  5,5\right)
$-division ring $\mathcal{D}_{M}^{\left(  5,5\right)  }$ (\ref{cn-f55}) is a
\textsf{quer-nonsymmetric} division ring.
\end{example}

\begin{definition}
\label{cn-def-dpf}The \textit{polyadic zeroless direct product field}
$\widehat{\mathcal{F}}^{\prime\left(  m,n\right)  }=\left\langle R^{\prime
}\mid\mathbf{\nu}^{\prime\left(  m\right)  },\mathbf{\mu}^{\prime\left(
n\right)  }\right\rangle $ consists of (two) zeroless polyadic fields
$\widehat{\mathcal{F}}_{1}^{\left(  m,n\right)  }=\left\langle R_{1}\mid
\nu_{1}^{\left(  m\right)  },\mu_{1}^{\left(  n\right)  }\right\rangle $ and
$\widehat{\mathcal{F}}_{2}^{\left(  m,n\right)  }=\left\langle R_{2}\mid
\nu_{2}^{\left(  m\right)  },\mu_{2}^{\left(  n\right)  }\right\rangle $ of
\textsf{the same} arity shape, while the componentwise operations on the
doubles $\mathbf{X}\in R_{1}\times R_{2}$ in (\ref{cn-y1})--(\ref{cn-y2})
still hold valid, and $\left\langle R_{1}\mid\mu_{1}^{\left(  n\right)
}\right\rangle $, $\left\langle R_{2}\mid\mu_{2}^{\left(  n\right)
}\right\rangle $, $\left\langle R^{\prime}=\left\{  \mathbf{X}\right\}
\mid\mathbf{\mu}^{\prime\left(  n\right)  }\right\rangle $ are $n$-ary groups.
\end{definition}

Following \textbf{Definition} \ref{cn-def-fmn1}, we have

\begin{corollary}
If at least one of the constituent fields is semicommutative, and another one
is totally commutative, then the polyadic product will be a semicommutative
$\left(  m,n\right)  $-field.
\end{corollary}

The additive and multiplicative unary queroperations (\ref{cg-m1g}) and
(\ref{cn-v1}) for the direct product field $\widehat{\mathcal{F}}^{\left(
m,n\right)  }$ are defined componentwise on the doubles $\mathbf{X}$ as
follows%
\begin{align}
\mathbf{\tilde{\nu}}^{\prime\left(  1\right)  }\left[  \mathbf{X}\right]   &
=\left(
\begin{array}
[c]{c}%
\tilde{\nu}_{1}^{\left(  1\right)  }\left[  x_{1}\right] \\
\tilde{\nu}_{2}^{\left(  1\right)  }\left[  x_{2}\right]
\end{array}
\right)  ,\label{cn-nx1}\\
\mathbf{\bar{\mu}}^{\prime\left(  1\right)  }\left[  \mathbf{X}\right]   &
=\left(
\begin{array}
[c]{c}%
\bar{\mu}_{1}^{\left(  1\right)  }\left[  x_{1}\right] \\
\bar{\mu}_{2}^{\left(  1\right)  }\left[  x_{2}\right]
\end{array}
\right)  ,\ \ \ x_{1}\in R_{1},\ x_{2}\in R_{2}. \label{cn-nx2}%
\end{align}

\begin{definition}
\label{cn-def-qsym1}A polyadic direct product field $\widehat{\mathcal{F}%
}^{\prime\left(  m,n\right)  }=\left\langle R^{\prime}\mid\mathbf{\nu}%
^{\prime\left(  m\right)  },\mathbf{\tilde{\nu}}^{\prime\left(  1\right)
},\mathbf{\mu}^{\prime\left(  n\right)  },\mathbf{\bar{\mu}}^{\prime\left(
1\right)  }\right\rangle $ is called \textit{quer-symmetric}, if its unary
queroperations (\ref{cn-nx1})--(\ref{cn-nx2}) commute%
\begin{align}
\mathbf{\tilde{\nu}}^{\prime\left(  1\right)  }\circ\mathbf{\bar{\mu}}%
^{\prime\left(  1\right)  }  &  =\mathbf{\bar{\mu}}^{\prime\left(  1\right)
}\circ\mathbf{\tilde{\nu}}^{\prime\left(  1\right)  },\label{cn-qsym1}\\
\widetilde{\overline{\mathbf{X}}}  &  =\overline{\widetilde{\mathbf{X}}%
},\ \ \ \forall\mathbf{X}\in R^{\prime}, \label{cn-qsym2}%
\end{align}
in other case $\widehat{\mathcal{F}}^{\prime\left(  m,n\right)  }$ is called a
\textit{quer-nonsymmetric direct product} $\left(  m,n\right)  $%
-\textit{field}.
\end{definition}

\begin{example}
Consider two nonunital zeroless $\left(  3,3\right)  $-fields $\widehat
{\mathcal{F}}_{1,2}^{\left(  3,3\right)  }=\left\langle \left\{
ia_{1,2}/b_{1,2}\right\}  \mid\nu_{1,2}^{\left(  3\right)  },\mu
_{1,2}^{\left(  3\right)  },\tilde{\nu}_{1,2}^{\left(  1\right)  },\bar{\mu
}_{1,2}^{\left(  1\right)  }\right\rangle $, $\mathsf{i}^{2}=-1$,
$a_{1,2},b_{1,2}\in\mathbb{Z}^{odd}$, where ternary additions $\nu
_{1,2}^{\left(  3\right)  }$ and ternary multiplications $\mu_{1,2}^{\left(
3\right)  }$ are sum and product in $\mathbb{Z}^{odd}$, correspondingly, and
the unary additive and multiplicative queroperations are $\tilde{\nu}%
_{1,2}^{\left(  1\right)  }\left[  \mathsf{i}a_{1,2}/b_{1,2}\right]
=-\mathsf{i}a_{1,2}/b_{1,2}$ and $\bar{\mu}_{1,2}^{\left(  1\right)  }\left[
\mathsf{i}a_{1,2}/b_{1,2}\right]  =-\mathsf{i}b_{1,2}/a_{1,2}$ (see
\textit{Example} \ref{cf-ex-tern}). Using (\ref{cn-y1})--(\ref{cn-y2}) we
build the operations of the polyadic nonderived nonunital zeroless product
$\left(  3,3\right)  $-field $\widehat{\mathcal{F}}^{\prime\left(  3,3\right)
}=\widehat{\mathcal{F}}_{1}^{\left(  3,3\right)  }\times\widehat{\mathcal{F}%
}_{2}^{\left(  3,3\right)  }$ on the doubles $\mathbf{X}^{T}=\left(
\mathsf{i}a_{1}/b_{1},\mathsf{i}a_{2}/b_{2}\right)  $ as follows%
\begin{equation}
\mathbf{\nu}^{\prime\left(  3\right)  }\left[  \mathbf{X}^{\left(  1\right)
},\mathbf{X}^{\left(  2\right)  },\mathbf{X}^{\left(  3\right)  }\right]
=\left(
\begin{array}
[c]{c}%
\mathsf{i}\dfrac{a_{1}^{\left(  1\right)  }b_{1}^{\left(  2\right)  }%
b_{1}^{\left(  3\right)  }+b_{1}^{\left(  1\right)  }a_{1}^{\left(  2\right)
}b_{1}^{\left(  3\right)  }+b_{1}^{\left(  1\right)  }b_{1}^{\left(  2\right)
}a_{1}^{\left(  3\right)  }}{b_{1}^{\left(  1\right)  }b_{1}^{\left(
2\right)  }b_{1}^{\left(  3\right)  }}\\[10pt]%
\mathsf{i}\dfrac{a_{2}^{\left(  1\right)  }b_{2}^{\left(  2\right)  }%
b_{2}^{\left(  3\right)  }+b_{2}^{\left(  1\right)  }a_{2}^{\left(  2\right)
}b_{2}^{\left(  3\right)  }+b_{2}^{\left(  1\right)  }b_{2}^{\left(  2\right)
}a_{2}^{\left(  3\right)  }}{b_{2}^{\left(  1\right)  }b_{2}^{\left(
2\right)  }b_{2}^{\left(  3\right)  }}%
\end{array}
\right)  ,
\end{equation}%
\begin{equation}
\mathbf{\mu}^{\prime\left(  3\right)  }\left[  \mathbf{X}^{\left(  1\right)
},\mathbf{X}^{\left(  2\right)  },\mathbf{X}^{\left(  3\right)  }\right]
=\left(
\begin{array}
[c]{c}%
-\mathsf{i}\dfrac{a_{1}^{\left(  1\right)  }a_{1}^{\left(  2\right)  }%
a_{1}^{\left(  3\right)  }}{b_{1}^{\left(  1\right)  }b_{1}^{\left(  2\right)
}b_{1}^{\left(  3\right)  }}\\[10pt]%
-\mathsf{i}\dfrac{a_{2}^{\left(  1\right)  }a_{2}^{\left(  2\right)  }%
a_{2}^{\left(  3\right)  }}{b_{2}^{\left(  1\right)  }b_{2}^{\left(  2\right)
}b_{2}^{\left(  3\right)  }}%
\end{array}
\right)  ,\ \ \ \ \ a_{i}^{\left(  j\right)  },b_{i}^{\left(  j\right)  }%
\in\mathbb{Z}^{odd},
\end{equation}
and the unary additive and multiplicative queroperations (\ref{cn-nx1}%
)--(\ref{cn-nx2}) of the direct product $\widehat{\mathcal{F}}^{\prime\left(
3,3\right)  }$ are%
\begin{align}
\mathbf{\tilde{\nu}}^{\prime\left(  1\right)  }\left[  \mathbf{X}\right]   &
=\left(
\begin{array}
[c]{c}%
-\mathsf{i}\dfrac{a_{1}}{b_{1}}\\[10pt]%
-\mathsf{i}\dfrac{a_{2}}{b_{2}}%
\end{array}
\right)  ,\\
\mathbf{\bar{\mu}}^{\prime\left(  1\right)  }\left[  \mathbf{X}\right]   &
=\left(
\begin{array}
[c]{c}%
-\mathsf{i}\dfrac{b_{1}}{a_{1}}\\[10pt]%
-\mathsf{i}\dfrac{b_{2}}{a_{2}}%
\end{array}
\right)  ,\ \ \ \ \ \ a_{i},b_{i}\in\mathbb{Z}^{odd}.
\end{align}

Therefore, both $\left\langle \left\{  \mathbf{X}\right\}  \mid\mathbf{\nu
}^{\prime\left(  3\right)  },\mathbf{\tilde{\nu}}^{\prime\left(  1\right)
}\right\rangle $ and $\left\langle \left\{  \mathbf{X}\right\}  \mid
\mathbf{\mu}^{\prime\left(  3\right)  },\mathbf{\bar{\mu}}^{\prime\left(
1\right)  }\right\rangle $ are commutative ternary groups, which means that
the polyadic direct product $\widehat{\mathcal{F}}^{\prime\left(  3,3\right)
}=\widehat{\mathcal{F}}_{1}^{\left(  3,3\right)  }\times\widehat{\mathcal{F}%
}_{2}^{\left(  3,3\right)  }$ is the \textsf{nonunital zeroless} polyadic
field. Moreover, $\widehat{\mathcal{F}}^{\prime\left(  3,3\right)  }$ is
quer-symmetric, because (\ref{cn-qsym1})--(\ref{cn-qsym2}) hold valid%
\begin{equation}
\mathbf{\bar{\mu}}^{\prime\left(  1\right)  }\circ\mathbf{\tilde{\nu}}%
^{\prime\left(  1\right)  }\left[  \mathbf{X}\right]  =\mathbf{\tilde{\nu}%
}^{\prime\left(  1\right)  }\circ\mathbf{\bar{\mu}}^{\prime\left(  1\right)
}\left[  \mathbf{X}\right]  =\left(
\begin{array}
[c]{c}%
\mathsf{i}\dfrac{b_{1}}{a_{1}}\\[10pt]%
\mathsf{i}\dfrac{b_{2}}{a_{2}}%
\end{array}
\right)  ,\ \ \ \ \ \ a_{i},b_{i}\in\mathbb{Z}^{odd}.
\end{equation}

\end{example}

\begin{example}
Let us consider the polyadic direct product of two zeroless fields, one of
them the semicommutative $\left(  5,5\right)  $-field $\widehat{\mathcal{F}%
}_{1}^{\left(  5,5\right)  }=\mathcal{F}_{M}^{\left(  5,5\right)  }$ from
(\ref{cn-f55}), and the other one the nonderived nonunital zeroless $\left(
5,5\right)  $-field of fractions $\widehat{\mathcal{F}}_{2}^{\left(
5,5\right)  }=\left\langle \left\{  \sqrt{\mathsf{i}}\frac{4r+1}%
{4s+1}\right\}  \mid\nu_{2}^{\left(  5\right)  },\mu_{2}^{\left(  5\right)
}\right\rangle $, $r,s\in\mathbb{Z}$, $\mathsf{i}^{2}=-1$. The double is
$\mathbf{X}^{T}=\left(  \sqrt{\mathsf{i}}\frac{4r+1}{4s+1},M\right)  $, where
$M$ is in (\ref{cn-m4}). The polyadic nonunital zeroless direct product field
$\widehat{\mathcal{F}}^{\prime\left(  5,5\right)  }=\widehat{\mathcal{F}}%
_{1}^{\left(  5,5\right)  }\times\widehat{\mathcal{F}}_{2}^{\left(
5,5\right)  }$ is nonderived and semicommutative, and is defined by
$\widehat{\mathcal{F}}^{\left(  5,5\right)  }=\left\langle \mathbf{X}%
\mid\mathbf{\nu}^{\prime\left(  5\right)  },\mathbf{\mu}^{\prime\left(
5\right)  },\mathbf{\tilde{\nu}}^{\prime\left(  1\right)  },\mathbf{\bar{\mu}%
}^{\left(  1\right)  }\right\rangle $, where its addition and multiplication
are%
\begin{equation}
\mathbf{\nu}^{\prime\left(  5\right)  }\left[  \mathbf{X}^{\left(  1\right)
},\mathbf{X}^{\left(  2\right)  },\mathbf{X}^{\left(  3\right)  }%
,\mathbf{X}^{\left(  4\right)  },\mathbf{X}^{\left(  5\right)  }\right]
=\left(
\begin{array}
[c]{c}%
\sqrt{\mathsf{i}}\dfrac{4R_{\nu}+1}{4S_{\nu}+1}\\[10pt]%
\left(
\begin{array}
[c]{cccc}%
0 & \dfrac{4K_{\nu,1}+3}{4L_{\nu,1}+3} & 0 & 0\\
0 & 0 & \dfrac{4K_{\nu,2}+3}{4L_{\nu,2}+3} & 0\\
0 & 0 & 0 & \dfrac{4K_{\nu,3}+3}{4L_{\nu,3}+3}\\
\dfrac{4K_{\nu,4}+3}{4L_{\nu,4}+3} & 0 & 0 & 0
\end{array}
\right)
\end{array}
\right)  ,
\end{equation}%
\begin{equation}
\mathbf{\mu}^{\prime\left(  5\right)  }\left[  \mathbf{X}^{\left(  1\right)
},\mathbf{X}^{\left(  2\right)  },\mathbf{X}^{\left(  3\right)  }%
,\mathbf{X}^{\left(  4\right)  },\mathbf{X}^{\left(  5\right)  }\right]
=\left(
\begin{array}
[c]{c}%
\sqrt{\mathsf{i}}\dfrac{4R_{\mu}+1}{4S_{\mu}+1}\\[10pt]%
\left(
\begin{array}
[c]{cccc}%
0 & \dfrac{4K_{\mu,1}+3}{4L_{\mu,1}+3} & 0 & 0\\
0 & 0 & \dfrac{4K_{\mu,2}+3}{4L_{\mu,2}+3} & 0\\
0 & 0 & 0 & \dfrac{4K_{\mu,3}+3}{4L_{\mu,3}+3}\\
\dfrac{4K_{\mu,4}+3}{4L_{\mu,4}+3} & 0 & 0 & 0
\end{array}
\right)
\end{array}
\right)  ,
\end{equation}
where $R_{\nu,\mu},S_{\nu,\mu}\in\mathbb{Z}$ are cumbersome integer functions
of $r^{\left(  i\right)  },s^{\left(  i\right)  }\in\mathbb{Z}$,
$i=1,\ldots,5$, and $K_{\nu,i},K_{\mu,i},L_{\nu,i},L_{\mu,i}\in\mathbb{Z}$ are
cumbersome integer functions of $k_{j}^{\left(  i\right)  },l_{j}^{\left(
i\right)  }\in\mathbb{Z}$, $j=1,\ldots,4$, $i=1,\ldots,5$ (see (\ref{cn-m4})).
The unary queroperations (\ref{cn-nx1})--(\ref{cn-nx2}) of the direct product
$\widehat{\mathcal{F}}^{\left(  5,5\right)  }$ are%
\begin{align}
\mathbf{\tilde{\nu}}^{\prime\left(  1\right)  }\left[  \mathbf{X}\right]   &
=\left(
\begin{array}
[c]{c}%
-3\sqrt{\mathsf{i}}\dfrac{4r+1}{4s+1}\\[10pt]%
-3M
\end{array}
\right)  ,\\
\mathbf{\bar{\mu}}^{\prime\left(  1\right)  }\left[  \mathbf{X}\right]   &
=\left(
\begin{array}
[c]{c}%
-\sqrt{\mathsf{i}}\left(  \dfrac{4s+1}{4r+1}\right)  ^{3}\\[10pt]%
\dfrac{4l_{1}+3}{4k_{1}+3}\dfrac{4l_{2}+3}{4k_{2}+3}\dfrac{4l_{3}+3}{4k_{3}%
+3}\dfrac{4l_{4}+3}{4k_{4}+3}M
\end{array}
\right)  ,\ \ \ \ \ r,s,k_{i},l_{i}\in\mathbb{Z},
\end{align}
where $M$ is in (\ref{cn-m4}). Therefore, $\left\langle \left\{
\mathbf{X}\right\}  \mid\mathbf{\nu}^{\prime\left(  5\right)  },\mathbf{\tilde
{\nu}}^{\prime\left(  1\right)  }\right\rangle $ is a commutative $5$-ary
group, and $\left\langle \left\{  \mathbf{X}\right\}  \mid\mathbf{\mu}%
^{\prime\left(  5\right)  },\mathbf{\bar{\mu}}^{\prime\left(  1\right)
}\right\rangle $ is a semicommutative $5$-ary group, which means that the
polyadic direct product $\widehat{\mathcal{F}}^{\prime\left(  5,5\right)
}=\widehat{\mathcal{F}}_{1}^{\left(  5,5\right)  }\times\widehat{\mathcal{F}%
}_{2}^{\left(  5,5\right)  }$ is the \textsf{nonunital zeroless} polyadic
semicommutative $\left(  5,5\right)  $-field. Using (\ref{cf-m81}) we obtain%
\begin{equation}
\mathbf{\tilde{\nu}}^{\prime\left(  1\right)  }\mathbf{\bar{\mu}}%
^{\prime\left(  1\right)  }\left[  \mathbf{X}\right]  =81\mathbf{\bar{\mu}%
}^{\prime\left(  1\right)  }\mathbf{\tilde{\nu}}^{\prime\left(  1\right)
}\left[  \mathbf{X}\right]  ,
\end{equation}
and therefore the direct product $\left(  5,5\right)  $-field $\widehat
{\mathcal{F}}^{\prime\left(  5,5\right)  }$ is quer-nonsymmetric (see
(\ref{cn-n1m})).
\end{example}

Thus, we arrive at

\begin{theorem}
The category of \textsf{zeroless }polyadic fields $\mathtt{zlessPolField}$ can
exist (having the class of all zeroless polyadic fields for objects and field
homomorphisms for morphisms) and can be well-defined, because it has a product
as the polyadic field product.
\end{theorem}

Further analysis of the direct product constructions introduced here and their
examples for polyadic rings and fields would be interesting to provide in
detail, which can also lead to new kinds of categories.

\bigskip

\textbf{Acknowledgement}. The author is deeply grateful to Vladimir Akulov,
Mike Hewitt, Vladimir Tkach, Raimund Vogl and Wend Werner for numerous
fruitful discussions, important help and valuable support.

\newpage

\pagestyle{emptyf}

\end{document}